\documentclass[reqno,12pt]{amsart}
\usepackage{amssymb,amsmath,amstext,amsthm,amsfonts}
\usepackage[ansinew]{inputenc}
\usepackage{a4wide}

\title[Stochastic stability of
diffeomorphisms]{Stochastic stability of non-uniformly\\
  hyperbolic diffeomorphisms}

\date{\today}

\begin{thanks} {J.F.A. and V.A. were partially supported by
    FCT through CMUP, POCI/MAT/61237/2004, and grants of
    FCT.  C.V. was partially supported by Fondecyt 11060497
    and DGIP-UCN 2005 10301201 and V.A. was also partially
    supported by CNPq.  Part of this work was done while
    V.A. enjoyed a leave at IMPA and IM-UFRJ and C.V enjoyed
    a leave at CMUP.}
\end{thanks}

\author{José F. Alves}
\address{José F. Alves, Departamento de Matemática Pura,
  Rua do Campo Alegre 687, 4169-007 Porto, Portugal}
\email{jfalves@fc.up.pt} \urladdr{http://www.fc.up.pt/cmup/jfalves}

\author{Vítor Ara\'ujo}
\address{Vítor Ara\'ujo, Instituto de Matemática,
  Universidade Federal do Rio de Janeiro, C. P. 68.530,
  21.945-970 Rio de Janeiro, RJ-Brazil \emph{and}
Centro de Matem\'atica da Universidade do Porto,
  Rua do Campo Alegre 687, 4169-007 Porto, Portugal}
\email{vitor.araujo@im.ufrj.br \text{and} vdaraujo@fc.up.pt}

\author{Carlos H. Vásquez}
\address{Carlos H. Vásquez,
Departamento de Matemática,
Universidad Catolica del Norte,
Av. Angamos 0610,
Antofagasta-Chile
}
\email{cavasquez@ucn.cl}

\subjclass{Primary: 37D25. Secondary: 37C40, 37H15.}

\keywords{Dominated splitting, non-uniform hyperbolicity, SRB
measure, random perturbation, stochastic stability}

\begin{document}

\begin{abstract} We prove that the statistical properties of
  random perturbations of a
  diffeomorphism with dominated splitting having mostly
  contracting center-stable direction and
  non-uniformly expanding center-unstable direction
  are described by a finite number of
  stationary measures. We also give necessary and sufficient
  conditions for the stochastic stability of such dynamical
  systems. We show that a certain $C^2$-open class of
  non-uniformly hyperbolic diffeomorphisms introduced by Alves,
  Bonatti and Viana in~\cite{ABV}
  are stochastically stable.
\end{abstract}

\newcommand{\mcup}{\mbox{$\bigcup$}}
\newcommand{\mcap}{\mbox{$\bigcap$}}

\def \RR {{\mathbb R}}
\def \ZZ {{\mathbb Z}}
\def \NN {{\mathbb N}}
\def \PP {{\mathbb P}}
\def \TT {{\mathbb T}}
\def \CC {{\mathbb C}}
\def \DD {{\mathbb D}}
\def \PP {{\mathbb P}}

 \def \neg {\!\!\!\!\!\!\!\!}
 \def \lra {\longrightarrow }
 \def \ra {\rightarrow }
 \def \wh {\widehat }
 \def \un {\underline }
 \def \dist {\operatorname{dist}}
 \def \dim {\operatorname{dim}}
 \def \ov {\overline}
 \def \m {\operatorname{m}}
 \def \inter{\mbox{int}\: }
 \def \supp {\operatorname{supp}\,}
 \def \wlim {\mbox{$w^*$-}\lim_{n\ra\infty}\, }
 \def \distp {\mbox{d}_\PP }
 \def \Leb {\operatorname{Leb}}
 \def \jac {\operatorname{Jac}}

\def \al {\alpha } \def \be {\beta }    \def \de {\delta }
\def \ga {\gamma } \def \ep {\varepsilon } \def \vfi {\varphi}
\def \th {\theta } \def \si {\sigma }   \def \om {\omega}
\def \la {\lambda} \def \ze {\zeta}

\def \La {\Lambda} \def \De {\Delta}

 \def \F {\mathcal{F}}
 \def \B {\mathcal{B}}
 \def \R {\mathcal{R}}
 \def \M {\mathcal{M}}
 \def \N {\mathcal{N}}
 \def \Q {\mathcal{Q}}
 \def \P {\mathcal{P}}
 \def \C {\mathcal{C}}
 \def \H {\mathcal{H}}
 \def \L {\mathcal{L}}
 \def \V {\mathcal{V}}
 \def \U {\mathcal{U}}
 \def \D {\mathcal{D}}
 \def \K {\mathcal{K}}

\newcommand{\dem}{\begin{proof}}
\newcommand{\cqd}{\end{proof}}

\newcommand{\qand}{\quad\text{and}\quad}

\newtheorem{maintheorem}{Theorem}
\renewcommand{\themaintheorem}{\Alph{maintheorem}}
\newcommand{\cmt}{\begin{maintheorem}}
\newcommand{\fmt}{\end{maintheorem}}

\newtheorem{maincorollary}[maintheorem]{Corollary}
\renewcommand{\themaintheorem}{\Alph{maintheorem}}
\newcommand{\cmc}{\begin{maincorollary}}
\newcommand{\fmc}{\end{maincorollary}}

\newtheorem{T}{Theorem}[section]
\newcommand{\ct}{\begin{T}}
\newcommand{\ft}{\end{T}}

\newtheorem{Corollary}[T]{Corollary}
\newcommand{\cco}{\begin{Corollary}}
\newcommand{\fco}{\end{Corollary}}

\newtheorem{Proposition}[T]{Proposition}
\newcommand{\cpr}{\begin{Proposition}}
\newcommand{\fpr}{\end{Proposition}}

\newtheorem{Lemma}[T]{Lemma}
\newcommand{\cle}{\begin{Lemma}}
\newcommand{\fle}{\end{Lemma}}

\theoremstyle{remark}

\newtheorem{Remark}[T]{Remark}
\newcommand{\cre}{\begin{Remark}}
\newcommand{\fre}{\end{Remark}}

\newtheorem{Definition}[T]{Definition}
\newcommand{\cd}{\begin{Definition}}
\newcommand{\fd}{\end{Definition}}

\newtheorem{Example}[T]{Example}
\newcommand{\cex}{\begin{Example}}
\newcommand{\fex}{\end{Example}}

\maketitle

\setcounter{tocdepth}{1}

\tableofcontents

\section{Introduction}
\label{sec:introduction}

Dynamical Systems Theory focuses mainly on two problems: what is the
behavior of a given system, and how this behavior changes under
small modifications of the law that governs the system. Properties
which are shared by the original system and its perturbations are
called ``stable''. This work is about the stability of certain
dynamical systems from a probabilistic point of view.

Much of the recent progress in Dynamics  is a consequence of  a
probabilistic approach to the understanding of complicated dynamical
systems, where one focuses on the statistical properties of
``typical orbits'', in the sense of large volume in the ambient
space. We deal here with diffeomorphisms $f:M\to M$ on compact
manifolds. The most basic statistical data are the time averages

$$
\lim_{n\to\infty} \frac 1n \sum_{j=0}^{n-1} \delta_{f^j(z)} $$

\noindent where $\delta_w$ represents the Dirac measure at a point
$w$.  The Ergodic Theorem asserts that time averages do exist at
almost every point for any invariant probability.  Moreover, if the
measure is ergodic then the time average coincides with the space
average, that is, the invariant probability measure itself. However,
many invariant measures are singular with respect to volume in
general, and so the Ergodic Theorem is not enough to understand the
behavior of positive volume (Lebesgue measure) sets of orbits.

An \emph{SRB measure} is an invariant probability measure for which
time averages exist and coincide with the space average, for a set
of initial conditions with positive Lebesgue measure. This set is
called the \emph{basin} of the measure. Sinai, Ruelle and Bowen
 introduced this notion  more than thirty years ago,
and proved that for uniformly hyperbolic (Axiom~A) diffeomorphisms
and flows time averages exist for Lebesgue almost every point and
coincide with one of finitely many SRB measures;
see~\cite{S72,BR,Bow,R76}.

The problem of existence and finiteness of SRB measures and their
stability under perturbations of the system, beyond the Axiom A
setting, remains a main goal of Dynamics. The construction of the so
called Gibbs $u$-states, by Pesin and Sinai in \cite{PS82} was the
beginning of the extension of the Sinai, Ruelle and Bowen ideas to
partially hyperbolic systems, a fruitful generalization of the
notion of uniform hyperbolicity, which more recently was shown to
encompass Lorenz-like and singular-hyperbolic flows
\cite{Tu,MPP99,MPP} and to be a consequence of robust transitivity
\cite{BDP,MPP}. We refer the reader to \cite{BDV,Vi98,Yo02} for
surveys on much of the progress obtained so far.

The papers of Alves, Bonatti and Viana ~\cite{ABV,BV},  and
Dolgopyat ~\cite{Do01} are of special interest to us here since they
prove existence and finiteness of SRB measures for partially
hyperbolic diffeomorphisms, under the assumption that the central
direction is either ``mostly contracting'' \cite{BV, Do01} or
``non-uniformly expanding'' \cite{ABV}. Indeed, we are going to
prove that a large class of systems in this setting, having weak
hyperbolicity properties, namely an invariant dominated splitting
with weak expansion and contraction, are stochastically stable. This
kind of diffeomorphisms has been intensely studied in  recent years.
See ~\cite{BDV,Vi98} and references therein for a global view of the
state of the art in this setting.

It was conjectured by Palis~\cite{Pa00} that every dynamical system
can be approximated by another one having only finitely many
SRB~measures, whose basins cover Lebesgue almost every point of the
ambient manifold. Moreover, these SRB measures should be stable
under perturbations of the system, in a stochastic sense.

Stochastic stability means that time averages do not change much
under small random perturbations -- the precise definition will be
given in the next section. This concept was introduced by Kolmogorov
and Sinai and much developed by the pioneering work of Kifer; see
~\cite{Ki86a,Ki88} and references therein. The reader can also see
~\cite{Liu2} for a recent survey on random dynamical systems.

Stochastic stability is well-known for uniformly expanding maps and
for uniformly hyperbolic diffeomorphisms~\cite{Ki86a,Yo,Vi}
restricted to the basin of their attractors. More recently,
\cite{AA} proved that the same is true for an open class of
non-uniformly expanding maps.

In this paper we prove stochastic stability for a large class of
diffeomorphisms admitting an invariant dominated splitting with
non-uniformly expanding center-unstable direction and mostly
contracting center-stable direction. Precise definitions and
statements of our results follow.



\subsection{Non-uniformly hyperbolic diffeomorphisms}
\label{sec:part-hyperb-diff}

Let $f:M\ra M$ be a smooth map defined on a compact Riemannian
manifold $M$.  We write $\|\cdot\|$ for the induced norm on $TM$
and fix some normalized Riemannian volume form $m$ on $M$ that we
call {\em Lebesgue measure}.  As explained in the Introduction, we
say that an $f$-invariant Borel probability measure $\mu$ on $M$
is an \emph{SRB measure} if for a positive Lebesgue measure set of
points $x\in M$
\begin{equation}\label{average} \lim_{n\ra +\infty}
\frac{1}{n}\sum_{j=0}^{n-1}\vfi\big(f^n x\big) =\int \vfi\, d\mu,
\end{equation} for every continuous map $\vfi:M\to \RR$. We
define the {\em basin} $B(\mu)$ of $\mu$ as the set of those
points $x$ in $M$ for which (\ref{average}) holds for all
continuous $\vfi$.


Let $f:M\to M$ be a diffeomorphism for which there is a strictly
forward $f$-invariant open set $U\subset M$, meaning that
$f(\ov{U})\subset U$. We say that $\Lambda=\cap_{n\ge0}f^n(U)$, the
maximal invariant subset in $U$,  is an \emph{attractor}. We suppose
that there is a continuous splitting $T_U M=E^{cs}\oplus E^{cu}$ of
the tangent bundle over $U$ which is $Df$-invariant, that is
$Df(E^{cs}_x)=E^{cs}_{fx}$ and $Df^{-1}(E^{cu}_{fx})=E^{cu}_x $ for
every $x\in U$. The bundles $E^{cs}$ and $E^{cu}$ will be called
\emph{center-stable} and \emph{center-unstable} and have dimensions
$s$ and $u$, respectively, with $u,s\ge1$ and $s+u=\dim (M)$. We
will assume several conditions on the splitting of $T_UM$:
\begin{enumerate}
\item[(a)] \emph{Dominated decomposition:} there
exists a constant $0<\lambda<1$ for which
$$
\| Df\vert E^{cs}_x
\|\cdot\| Df^{-1}\vert E^{cu}_{fx} \| \le\lambda \quad\mbox{for
all $x\in  U$.}
$$
\item[(b)] \emph{Nonuniform expansion along the
central-unstable direction:} there is $c_u>0$ such that for
Lebesgue almost all $x\in U$
$$
  \limsup_{n\to+\infty}\frac1n\sum_{j=0}^{n-1} \log\| Df^{-1} \vert
  E^{cu}_{f^jx} \| \le -c_u.
$$
\end{enumerate}

Theorem 6.3 of \cite{ABV} shows that a diffeomorphism $f$
with a dominated splitting non-uniformly expanding along the
center-unstable direction has some ergodic \emph{Gibbs
  $cu$-state} $\mu$ supported in the attractor $\La$, that
is, $\mu$ is an invariant probability measure whose $u$
larger Lyapunov exponents are positive and whose conditional
measures along the corresponding local unstable manifolds
are almost everywhere absolutely continuous with respect to
Lebesgue measure on such manifolds. This notion is a
non-uniform version of the Gibbs $u$-states introduced by
Pesin and Sinai ~\cite{PS}. Moreover if the derivative of
$f$ uniformly contracting along the sub-bundle $E^{cs}$,
then $\mu$ is an SRB measure. This is a well known
consequence of the absolute continuity of the conditional
measures of $\mu$ and the absolute continuity of the stable
lamination, see e.g. \cite{Pe76,BDV}: the union of the
stable manifolds through the points whose time averages are
given by $\mu$ is a positive Lebesgue measure subset of the
basin of $\mu$.

As shown in \cite{BV, Do01} there are some cases where non-uniform
contraction along $E^{cs}$ is enough for ensuring the
existence of SRB measures for $f$. A sufficient condition which we
assume here is that on every disk at any local unstable manifold
there exists a positive Lebesgue measure subset of points of that
disk having negative Lyapunov exponents along the central-stable
direction:


\begin{enumerate}
\item[(c)] \emph{Mostly contracting
central-stable direction:} every disk $D$ contained in an
unstable local manifold satisfies
\begin{equation*}
\limsup_{n\to+\infty}\frac{1}{n}\log\|Df^n|E^{cs}_x\|<0
\end{equation*}
for a positive Lebesgue measure set of points $x\in D$.
\end{enumerate}

Conditions (a)-(c) above are the main conditions we assume
on the $C^2$ diffeomorphism~$f$.  Condition (a) is a
classical one and was already considered in \cite{HPS} as a
weakening of the uniform hyperbolicity conditions. Condition
(b) was first considered in \cite{ABV} and together with
condition (a) it ensures the existence of Gibbs $cu$-states.
Condition (c) was considered in \cite{BV, Do01} as a means
to ensure that Gibbs $cu$-states are SRB measures.  These
conditions are satisfied by $C^2$ open classes of
diffeomorphisms which are not uniformly hyperbolic systems,
as shown in the last section.



\subsection{Statement of results}
\label{sec:random-perturbations}

In this work we are interested in studying random
perturbations of a $C^2$ diffeomorphism $f\colon M\to M$
satisfying conditions (a)-(c) from
Subsection~\ref{sec:part-hyperb-diff}. For that, we take a
continuous map
 \[
 \begin{array}{rccl}
 \Phi:& T &\longrightarrow&  \mbox{Diff}^2(M,M)\\
 & t &\longmapsto & f_t
 \end{array}
 \]
from a metric space $T$ into  the space of $C^2$ diffeomorphisms
of  $M$, with $f=f_{t^*}$ for some fixed $t^*\in T$. Given $x\in
M$, a \emph{random orbit} of $x$ will be a sequence $(f_{\un t}^n
x)_{n\ge1}$ where $\un t$ denotes an element
$(t_1,t_2,t_3,\ldots)$ in the product space $T^{\NN}$ and
 $$
 f^n_{\un t}= f_{t_n}\circ \cdots \circ
 f_{t_1}\quad \mbox{for}\quad n \ge1.
 $$
 We also take a family $(\th_\ep)_{\ep>0}$ of probability
 measures on $T$ such that their supports $\supp(\th_\ep)$
 form a nested family of connected compact sets and
 $\supp(\th_\ep)\rightarrow \{t^*\}$ when $ \ep\to 0$.  We
 assume some quite general non-degeneracy conditions on the
 map $\Phi$ and the measures $\th_\ep$ (see the beginning of
 Section~\ref{sec:stochastic-stability}) and refer to
 $\{\Phi,(\th_\ep)_{\ep>0}\}$ as a {\em random perturbation}
 of $f$. It is known \cite{Ar1} that  every smooth map
$f:M\ra M$ of a compact manifold  always admits a random
perturbation satisfying these non-degeneracy conditions.

In the setting of random perturbations of a map, a Borel
probability measure $\mu^\ep$ on $M$ is said to be a
\emph{random physical measure} if for a positive Lebesgue
measure set of points $x\in M$ one has
 \begin{equation}\label{eq:pertaverage}
 \lim_{n\ra +\infty}
 \frac{1}{n}\sum_{j=0}^{n-1}\vfi\big(f_{\un t}^j x\big)
 =\int \vfi \,d\mu^\ep,
 \end{equation}
 for all continuous $\vfi\colon M\rightarrow\RR$ and
 $\th^\NN_\ep$ almost every $\un t\in T^\NN.$ We denote the
 set of points $x\in M$ for which (\ref{eq:pertaverage})
 holds for every continuous $\varphi$ and for $\th^\NN_\ep$
 almost every $\un t\in T^\NN$ by $B(\mu^\ep)$, and say that
 it is the \emph{basin of $\mu^\ep$}.

 The map $f\colon M\ra M$ is said to be {\em stochastically
   stable} if the weak$^*$ accumulation points (when $\ep>0$
 goes to zero) of the physical measures of the random
 perturbation are convex linear combinations of the SRB
 measures of $f$ (for this notion of stochastic stability,
 see~\cite{Ar2}).

 In this setting it is known that there exists a finite
 number $l$ (bounded above by the number of $SRB$ measures
 of $f$) of random physical measures whose basins cover the
 ambient space. That is, almost all random orbits of every
 point of $M$ have time averages given by a random physical
 measure out of $\mu^\ep_1,\dots,\mu^\ep_l$. See
 Section~\ref{sec:stochastic-stability} for a detailed
 statement of this result and an outline of its proof.  This
 property is akin to the one obtained in~\cite{ABV} for
 points in the topological basin of attractors with a
 dominated splitting having non-uniform expansion along the
 central-unstable direction and uniform contraction along
 the central-stable direction (usually referred to as
 \emph{partially} hyperbolic attractors, see
 e.g.~\cite{BDV}).

Now we present a notion that will play an important role in
the study of the stochastic stability of a diffeomorphism
with a dominated splitting. We say that $f$ is {\em
  non-uniformly expanding along the center-unstable
  direction for random orbits} if there is $c>0$ such that
\begin{equation}
  \label{eq:rNUE}
\limsup_{n\to+\infty}\frac1n\sum_{j=0}^{n-1} \log\| Df^{-1}
\vert E^{cu}_{f_{\un t}^jx} \| \le -c
\end{equation}
for Lebesgue almost every $x\in U$ and $\th_\ep^\NN$
almost every $\un t\in T^\NN$, at least for small $\ep>0$.

The main result below gives a characterization of the
stochastically stable non-uniformly hyperbolic
diffeomorphisms. To the best of our knowledge, currently
this is the only result pointing in the direction of general
stochastic stability for diffeomorphisms with a dominated
splitting.

\begin{maintheorem}\label{te:mainB}
  Let $f:M\to M$ be a $C^2$ diffeomorphism admitting a
  strictly forward invariant open set $U$ with a
  $Df$-invariant dominated splitting $T_U M=E^{cs}\oplus
  E^{cu}$, where the center-stable direction is mostly
  contracting and the center-unstable direction is
  non-uniformly expanding.  Then $f$ is stochastically
  stable if, and only if, $f$ is non-uniformly expanding
  along the center-unstable direction for random orbits.
 \end{maintheorem}

 In Section~\ref{sec:stoch-stable-class} we exhibit open
 classes of  diffeomorphisms satisfying the
 conditions of Theorem~\ref{te:mainB} which are
 stochastically stable. The ones we have in mind were
 presented in \cite[Appendix A]{ABV}. These form an open
 class $\D$ of diffeomorphisms defined on the
 $d$-dimensional torus $M=\TT^d$, $d\ge2$, having the whole
 $M$ as an attractor, $TM=E^{cs}\oplus E^{cu}$ with mostly
 contracting center-stable direction and exhibiting
 non-uniform expansion along the $E^{cu}$ direction.
In
 Section~\ref{sec:stoch-stable-class}
 we present a complete description and we prove that
 the members of this class are stochastically stable by
 showing that they satisfy the condition in
 Theorem~\ref{te:mainB} above.  Further developments in
 \cite{T} show that these diffeomorphisms have a unique SRB
 measure, and in \cite{Vz} it was already proved that these
 maps are also statistically stable.

 The conditions needed to get
 Theorem~\ref{te:mainB} are general enough to enable us to
 obtain as corollaries of our method the stochastic
 stability of some other families of diffeomorphisms with
 dominated splittings.
Obviously uniform expansion along the center-unstable bundle and
mostly contractive center-stable bundle fit in our conditions.
This setting was studied in~\cite{BV, Do01}, where it was shown
that these conditions are enough to obtain $SRB$ measures for $f$.
Since uniform expansion along the center-unstable direction and
dominated splitting are robust in a whole $C^1$-neighborhood of
$f$, we get non-uniform expansion for random orbits along the
center-unstable direction for free.

\cmc\label{cor.cs+uu} Let $f:M\to M$ be a $C^2$ diffeomorphism
admitting a strictly forward invariant open set $U$ with a dominated
splitting $T_U M=E^{cs}\oplus E^{u}$, where  the center-stable
direction  is mostly contracting and the center-unstable direction
is uniformly expanding. Then $f$ is stochastically stable. \fmc

Another easy remark is that uniform contraction along the
center-stable bundle with nonuniform expansion along the
center-unstable direction, together with nonuniform expansion for
random orbits along this direction, also fit in our setting.

\cmc\label{cor.ss+cu} Let $f:M\to M$ be a $C^2$ diffeomorphism
admitting a strictly forward invariant open set $U$ with a dominated
splitting $T_U M=E^{s}\oplus E^{cu}$, where the center-stable
direction is uniformly contracting and  the center-unstable
direction is non-uniformly expanding. If $f$ is non-uniformly
expanding on random orbits along the center-unstable direction, then
$f$ is stochastically stable. \fmc

Indeed, if $Df$ acting on $E^{cs}=E^s$ over $U$ is uniformly
contracting, then \emph{every} point of any disk $D$ inside
an unstable local manifold $W^u_{\mathrm{loc}}(x)$ for
$x\in\Lambda$ has negative Lyapunov exponents along the
$E^s$ direction. This is a simple consequence of the
definition of attractor: the unstable manifold of every
$x\in\Lambda$ is contained in $\Lambda$, so $D\subset
W^u_{\mathrm{loc}}(x)\subset\Lambda\subset U$ and hence $Df$
contracts $E^s_y$ for every $y\in D$.

In each of these corollaries the non-uniform bundle may admit a
further dominated splitting into a uniformly behaved bundle and a
central one, where the non-uniform expanding or contracting
conditions will apply. This is the setting of \emph{partially
hyperbolic diffeomorphisms} (see~\cite{BDV} for definitions and
main features), where the tangent bundle splits into three
sub-bundles
\[
T_U M=E^s\oplus E^c \oplus E^u
\]
such that both $(E^s\oplus E^c)\oplus E^u$ and $E^s\oplus (E^c \oplus
E^u)$ are dominated splittings, $E^s$ is uniformly contracted and
$E^u$ is uniformly expanding. Finally, the uniformly hyperbolic case
(when $E^c=\{0\}$) is also clearly in our setting (this result was
first proved in~\cite{Ki86a} but our approach is closer to \cite{Yo}).

We note that every statement is also true if the region $U$
coincides with the entire manifold, as the examples in
Section~\ref{sec:stoch-stable-class} and in~\cite{BV,ABV}.
In addition, the $C^2$ smoothness is not strictly necessary
for our arguments: $C^{1+\alpha}$ for some $\alpha>0$ is all
that is really needed in our arguments and in our
references.

The paper is organized as follows. We start the proofs by
studying random perturbations in
Section~\ref{sec:stochastic-stability}, where the proof of
the result on finiteness of random physical measures
(Theorem~\ref{te:mainA}) is explained. Then we prove the
necessary condition of Theorem~\ref{te:mainB} in
Section~\ref{sec:nonun-expans-rand}, where we also outline
the main steps of the proof of sufficiency.  The proof of
the sufficient condition of Theorem~\ref{te:mainB} starts at
Section~\ref{sec:main-lemmas} where we state some
preliminary results about hyperbolic times and distortion
control, extends to Section~\ref{sec:centre-unst-cylind}
where we study physical measures under random perturbations
with non-uniform expansion, and is completed in
Section~\ref{sec:accum-cylind} where we consider the limit
when the noise level tends to zero.  Finally, in
Section~\ref{sec:stoch-stable-class} we describe the
construction of an open class of non-uniformly hyperbolic
diffeomorphisms fitting the conditions of our theorems.


\section{Random perturbations}
\label{sec:stochastic-stability}


To begin  our study of random perturbations
$\{\Phi,(\th_\ep)_{\ep>0}\}$ of a diffeomorphism $f$
admitting a partially hyperbolic set $ U$, we observe that
if we take $\ep>0$ small enough, then all random orbits
$(f^n_{\un t}x)_{n\ge1}$ for every $x\in U$ are contained
in a compact subset $K$ of $ U$ when $\un t
\in\supp(\th_\ep^\NN)$, just by continuity of $\Phi$ and
because $U$ is strictly forward invariant.

Moreover if we introduce the  skew-product map $F:
T^\NN\times M \to T^\NN\times M$ given by $(\un t,
z)\mapsto(\sigma(\un t),f_{t_1}(z))$,
 where $\sigma$ is the left shift on
sequences $\un t=(t_1,t_2,\dots)\in T^\NN$, then we have
that
\begin{equation}
  \label{eq:randattractor}
\hat{\La}_\ep=\bigcap_{n\ge0} F^n( \supp(\th_\ep^\NN)\times U ),
\quad\mbox{with}\quad \La_\ep=\pi(\hat{\La}_\ep)\subset K,
\end{equation}
where $\pi:T^\NN\times M \to M$ is the natural projection.
$ \hat{\La}_\ep$ is a forward $F$ invariant set, the
\emph{$\ep$-random attractor}.

As mentioned before, we will assume that the random
perturbations of the partially hyperbolic map $f$ satisfy
some {\em non-degeneracy conditions}: there is $0<\ep_0<1$
such that for every $0<\ep<\ep_0$ we may find
$n_0=n_0(\ep)\in\NN$ satisfying the following conditions
for all $x\in M$ and all $n\ge n_0$:

 \begin{enumerate}
 \item  there is $\xi=\xi(\ep)>0$ such
 that $ \left\{ f^n_{\un t}x \colon \un{t}\in
 (\supp\th_\ep)^{\NN} \right\}\supset B(f^n x,\xi)$;
 \item defining $f^n_\odot x:T^{\NN}\ra M$  as $f^n_\odot x(\un t)=f_{\un
t}^n(x)$, then $(f^n_\odot x)_*{\th}^{\NN}_\ep\ll \m$.
 \end{enumerate}
 Condition~1 says that the set of
perturbed iterates of any given point  covers a full
neighborhood of the unperturbed iterates after a threshold
for all sufficiently small noise levels. Condition~2 means
that each set of perturbation vectors having positive
${\th}_{\ep}^{\NN}$ measure must send each point $x\in M$
onto positive Lebesgue measure subsets of $M$, after a
certain number of iterates.

Examples 1 and 2 in \cite{Ar1} show that \emph{every smooth
  map $f:M\ra M$ of a compact manifold always has a random
  perturbation satisfying the non-degeneracy conditions 1
  and 2}, as long as we take $T=\RR^p$, $t^*=0$ and also
$\th_\ep$ equal to Lebesgue measure restricted to the ball of radius
$\ep$ around $0$ (normalized to become a probability measure), for a
sufficiently big number $p\in\NN$ of parameters. For manifolds whose
tangent bundle is trivial (parallelizable manifolds) the random
perturbations which consist in adding at each step a random noise to
the unperturbed dynamics clearly satisfy non-degeneracy conditions 1
and 2 for $n_0=1$. We observe that this is precisely the kind of
perturbations considered by Baladi, Benedicks, Viana and Young in
\cite{BaV,BeV,BeY} for quadratic and Hénon maps.

The attractor $\cap_{n\ge0} f^n(U)$ for $f$ will be denoted
by $\La$. The conditions above imply that for small $\ep>0$
every $f_t$ is $C^2$-close to $f\equiv f_{t^*}$. Then
for every neighborhood $V$ of $\La$ we have $\La_\ep\subset
V$ for all small enough $\ep$. Thus the
compact set $K$ containing $\La_\ep$ may be taken as a
neighborhood of $\La$. We assume this in the rest of the
paper (this is important in Subsection~\ref{sec:sufficient}).


\subsection{The number of random physical measures}
\label{sec:finite-physical}

In the setting of random perturbations of a map, we say that
a set $A\subset M$ is {\em forward invariant} if for some
given small $\ep>0$ and for all $t\in\supp(\th_\ep)$ we have $f_t
A\subset A$. A probability measure $\mu$ is said to be
\emph{stationary}, if for every continuous $\vfi:M\to\RR$ it
holds the following relation, similar to the non-random
setting of invariant measures:
 \begin{equation}
 \label{eq.stationary}
 \int\vfi\,d\mu=
 \int\int\vfi\big(f_t x\big)\,d\mu(x)\,d\th_\ep(t).
 \end{equation}

\cre\label{re.accinvariant}
 If $(\mu^\ep)_{\ep>0}$ is a
family of stationary measures having $\mu_0$ as a weak$^*$
accumulation point when $\ep$ goes to $0$, then it follows
from~(\ref{eq.stationary}) and the convergence of
$\supp(\th_\ep)$ to $\{t^*\}$  that $\mu_0$ must be
invariant by $f=f_{t^*}$.
 \fre
 \noindent Condition~(\ref{eq.stationary}) is equivalent to
 saying that
 $F_{*}(\th_\ep^\NN\times\mu)=\th_\ep^\NN\times\mu$ and it
 is easy to see that a stationary measure $\mu$ satisfies $$
 x\in\supp(\mu) \quad\Rightarrow\quad f_t (x)
 \in\supp(\mu)\mbox{ for all }t\in\supp(\th_\ep), $$ just by
 continuity of $\Phi$. This means that if $\mu$ is a
 stationary measure, then $\supp(\mu)$ is a forward
 invariant set. Then the interior of $\supp(\mu)$ is
 nonempty, by non-degeneracy condition 1, and forward
 invariant by the continuity of the maps $f_t$. Obviously a
 physical measure is stationary.  In our setting the random
 physical measures $\mu$ we deal with are such that
 $\th^\NN_\ep\times\mu$ is $F$ ergodic.  Moreover the
 non-degeneracy condition 2 ensures that each stationary
 probability measure $\mu$ is absolutely continuous, that is
 $\mu\ll m$.

\ct
  \label{te:mainA} Let $\{\Phi,(\th_\ep)_{\ep>0}\}$ be a
  random perturbation of a $C^2$ diffeomorphism $f$
  admitting a strictly forward invariant open set $U$ with a
  $Df$-invariant dominated splitting $T_U M=E^{cs}\oplus
  E^{cu}$, where the center-stable direction is mostly
  contracting and the center-unstable direction is
  non-uniformly expanding. Then there is $l\in\NN$ such that
  for small enough $\ep>0$ there exist random physical
  measures $\mu^\ep_1,\dots,\mu^\ep_l$ with pairwise
  disjoint supports with $\supp(\mu_i^\ep)\subset
  B(\mu_i^\ep)$ for $ i=1,\dots,l$, satisfying
\begin{enumerate} \item for each $x\in M$ and $\th_\ep^\NN$ almost
every $\un t\in T^\NN$ there is $i=i(x,\un t)\in\{1,\dots,l\}$ such
that
 $$
 \mu_i^\ep=\lim_{n\to+\infty} \frac1n \sum_{j=1}^{n-1}
 \de_{f_{\un t}^j x}
 ; 
 $$
\item for $1\le i\le l$ and any probability measure $\eta$
  with support contained in $B(\mu^\ep_i)$ we have
 $$
 \mu_i^\ep=\lim_{n\to+\infty}\frac1n\sum_{j=0}^{n-1} (f^j_{\un t})_*
 \eta\quad
 \mbox{
 for $\th_\ep^\NN$ almost every $\un t\in T^\NN$;}
 $$
\item the number of random physical measures is at most the number
  of SRB measures for the unperturbed map $f$.
 \end{enumerate}
\ft

The following  is a general result from \cite{Ar1} (see
also~\cite[Section 2]{BK} for an alternative approach) for
random perturbations of a diffeomorphism satisfying the
non-degeneracy conditions $1$ and $2$ above, which implies
the first item of Theorem~\ref{te:mainA}.

\ct\label{pr:finitemeasures} Assume that
$\{\Phi,(\th_\ep)_{\ep>0}\}$ is a random perturbation of a
$C^1$ diffeomorphism $f$ satisfying non-degeneracy
conditions $1$ and $2$. Then for each $\ep>0$ there are
random physical measures $\mu^\ep_1,\dots\mu^\ep_{l(\ep)}$
which are absolutely continuous, and for each $x\in M$ there
is a $\th_\ep^\NN$ mod $0$ partition
$T_1(x),\dots,T_{l(\ep)}(x)$ of $T^\NN$ such that for $1\le
i\le {l(\ep)}$
 $$
 \mu_i^\ep=\lim_{n\to+\infty} \frac1n \sum_{j=1}^{n-1}
 \de_{f_{\un t}^j x} \quad \mbox{ for}\quad \un t \in
 T_i(x).
 $$
 Moreover the supports of the physical measures are
 pairwise disjoint, have nonempty interior and are contained
 in their basins: $\supp(\mu_i^\ep)\subset B(\mu_i^\ep),
 i=1,\dots,l(\ep)$.
 \ft

 Now we prove the second item of Theorem~\ref{te:mainA}.
 First we observe that since the basin of each physical
 measure contains its support, then it also has nonempty
 interior. Let $\mu^\ep=\mu^\ep_i$ be a physical measure and
 take any probability measure $\eta$ with support contained
 in $B(\mu^\ep_i)$. Given any continuous function $\varphi\colon
 M\to\RR$, we have for each $x\in \supp(\eta)$ and $\th_\ep$
 almost every $\un t\in T^\NN$
 $$
 \lim_{n\ra +\infty}
 \frac{1}{n}\sum_{j=0}^{n-1}\vfi(f_{\un t}^j x)
 =\int \vfi \,d\mu^\ep
 $$
by definition of $B(\mu^\ep)$. Taking integrals over
$\supp(\eta)$ with respect to $\eta$ on both sides of the
equality, the Dominated Convergence Theorem gives
 $$
 \lim_{n\ra +\infty}
 \frac{1}{n}\sum_{j=0}^{n-1}\int\vfi(f_{\un t}^j
 x) \, d\eta(x)
 =\int \vfi \,d\mu^\ep.
 $$
(Recall that we are taking integrals over the support of
the probability measure $\eta$). Since
 $$
 \int\vfi\circ f_{\un t}^j\, d\eta= \int\vfi\, d(f_{\un t}^j)_*\eta
 $$
for every  integer $j\ge 0$, we have proved item 2 of
Theorem~\ref{te:mainA}.

\medskip

 To conclude the
proof of Theorem~\ref{te:mainA} all we need is to show that
$l(\ep)$ is constant for every sufficiently small $\ep>0$ and that
$l$ is at most the number of SRB measures of $f$. We start by
recalling that, by assumption, $\{\supp(\th_\ep)\}_{\ep>0}$ is a
nested family of sets. This implies that if $\mu^\ep$ is a
physical measure, then $\supp(\mu^\ep)$ is forward invariant with
respect to $f_t$ for all $t\in \supp(\th_{\ep'})$, whenever
$0<\ep'<\ep$. Since non-degeneracy conditions are satisfied in
$\supp(\mu^\epsilon)$ by the probability measure $\th_{\ep'}$,
then Theorem~\ref{pr:finitemeasures} ensures that there exists at
least one physical measure $\mu^{\ep'}$ with $
\supp(\mu^{\ep'})\subset \supp(\mu^{\ep})$ for $0<\ep'<\ep$. This
shows that the number of physical measures is a non-decreasing
integer function of $\ep$ when $\ep\to 0$. Hence if we show that
$l$ is bounded from above, we conclude that $l$ is constant for
all small enough $\ep>0$.

We observe that $\supp(\mu^\ep)$ is forward invariant under
$f=f_{t^*}$ and, moreover, conditions (a)-(c) hold on $f\mid
\supp(\mu^\ep)$ because they hold Lebesgue almost everywhere in
$U$ (by assumption) and $\supp(\mu^\ep)$ has nonempty interior.
Thus \cite[Theorem C]{Vz} guarantees the existence of at least one
SRB measure $\mu$ with $\supp(\mu)\subset\supp(\mu^\ep)$.

We have seen that each support of a physical measure
$\mu^\ep$ must contain at least the support of one SRB
measure for the unperturbed map $f$. Since the number of
SRB measures is finite we have $l\le p$, where $p$ is the
number of those measures. This concludes the proof of
Theorem~\ref{te:mainA}.


\section{Stochastic stability}
\label{sec:nonun-expans-rand}

In this section we start the proof of Theorem~\ref{te:mainB}. Let
$f:M\to M$ be a $C^2$ diffeomorphism admitting a strictly forward
invariant open set $U$ with a dominated splitting $T_U
M=E^{cs}\oplus E^{cu}$, where  the center-stable direction  is
mostly contracting and the center-unstable direction is
non-uniformly expanding. We first prove that non-uniform expansion
along the center-unstable direction for random orbits is a
necessary condition for stochastic stability.

\subsection{Stochastic stability implies non-uniform
hyperbolicity}
\label{sec:necessary}

Taking $\ep>0$ small we know from Theorem~\ref{te:mainA}
that there is a finite number of random physical measures
$\mu_1^\ep,\dots \mu_l^\ep$ and, for each $x\in U$, there is
a $\th_\ep^\NN$ mod $0$ partition $T_1(x),\dots,T_l(x)$ of
$T^\NN$ for which $$ \mu_i^\ep=\lim_{n\to+\infty} \frac1n
\sum_{j=1}^{n-1}\de_{f_{\un t}^j(x)} \quad \mbox{for each
}\un t \in T_i(x).  $$ Furthermore, since $\log\| Df^{-1}
\vert E^{cu}_{x} \|$ is a continuous map, then we have for
every $x\in U$ and $\th^\NN_\ep$ almost every $\un t\in
T^\NN$ $$ \lim_{n\ra \infty}\frac{1}{n}\sum_{j=0}^{n-1}
\log\| Df^{-1} \vert E^{cu}_{f_{\un t}^jx} \|= \int\log\|
Df^{-1} \vert E^{cu}_{x} \|d\mu_i^\ep, $$ for some physical
measure $\mu_i^\ep$ with $1\le i\le l$.  Hence, for proving
that $f$ is non-uniformly expanding along the
center-unstable direction for random orbits it suffices to
show that there is $c_0>0$ such that if $1\le i\le l$ then,
for small $\ep>0$, $$ \int\log\| Df^{-1} \vert E^{cu}_{x}
\|d\mu^\ep_i<-c_0.  $$

The next result is proved in \cite[Lemma 5.1]{AA} for endomorphisms,
but the proof still works in our case, since it only uses weak*
convergence and the definition of physical and SRB measures.

 \cle\label{l.aprox} Let $\vfi\colon M\ra \RR$ be
 continuous.  Given $\de>0$ there is $\ep_0>0$ such that if
 $\ep\leq\ep_0$, then
 $
 \big|\int \vfi d\mu^\ep - \int\vfi d\mu_\ep\big|<\de,
 $
 for some linear convex combination  $\mu_\ep$ of the
 physical measures of $f$.
 \fle

Therefore there are real numbers
$w_1(\ep),\dots ,w_p(\ep)\geq 0$ such that
$w_1(\ep)+\cdots+w_p(\ep)=1$ and
$\mu_\ep=w_1(\ep)\mu_1+\cdots+w_p(\ep)\mu_p$. Since $\mu_i$
is a physical measure for $1\leq i\leq p$, we have for Lebesgue
almost every $x\in B(\mu_i)$
\[
\int\log\| Df^{-1} \vert E^{cu}_{x} \|d\mu_i=\lim_{n\ra
+\infty}\frac{1}{n} \sum_{j=0}^{n-1}\log\| Df^{-1} \vert
E^{cu}_{f^j(x)} \|\leq -c<0.\]
Hence
\[
\int\log\| Df^{-1}
\vert E^{cu}_{x} \|d\mu_\ep\leq -c,
\quad\mbox{so}\quad
 \int\log\| Df^{-1}
\vert E^{cu}_{x} \|d\mu^\ep_i\leq -c/2,
\]
applying Lemma~\ref{l.aprox} to $\vfi=\log\| Df^{-1} \vert
E^{cu} \|$, obtaining that $f$ is non-uniformly expanding
along the center-unstable direction for random orbits, 
%
as we want.



\subsection{Non-uniform hyperbolicity implies stochastic
  stability}
\label{sec:sufficient}

Let us explain why in our setting non-uniform expansion along the
center-unstable direction  for random orbits is a sufficient
condition for stochastic stability.

In order to prove that $f=f_{t^*}$ is stochastically stable,
it suffices to show that every weak$^*$ accumulation point
$\mu$ of any family $(\mu^\ep)_{\ep>0}$, where each
$\mu^\ep$ is a random physical measure of level $\ep$, is a
\emph{Gibbs $cu$-state} for $f$, that is
\begin{itemize}
\item $\mu$ is absolutely continuous with respect to the
  Lebesgue measure along local center-unstable disks, and
\item the Lyapunov exponents of $\mu$-a.e. point along the
  tangent directions to these disks are all strictly
  positive.
\end{itemize}
This follows from a combination of results. The first one is
a consequence of the techniques developed in \cite{BV,ABV}
and is detailed in \cite[Theorem C]{Vz}.

\ct\label{t.ABV-prop6.3} Let $f$ be a $C^2$ diffeomorphism
having an attractor exhibiting a
dominated splitting which is non-uniformly expanding along
the center-unstable direction and mostly contracting along
the center-stable direction.

Then every ergodic Gibbs cu-state is an SRB measure for $f$.
Moreover there are finitely many SRB measures whose basins
cover a full Lebesgue measure subset of the topological
basin $U$ of the attractor.  \ft

The proof of the next result can be found in \cite[Corollary
4.1]{Vz}.

\ct\label{t.Vz-ThmA+CorD} Let $f$ be a $C^2$ diffeomorphism
having an attractor with dominated splitting and let $\mu$
be a Gibbs cu-state for $f$. Then every ergodic component of
$\mu$ is a Gibbs cu-state.

Assuming further that $f$ is non-uniformly expanding along
the center-unstable direction, then every ergodic SRB
measure is a Gibbs cu-state.
\ft

These results ensure that every Gibbs cu-state $\mu$ for $f$ has
finitely many ergodic components which are SRB measures for
$f$ and also that $\mu$ can be written as a linear convex
combination of these SRB measures.
Hence it is enough to prove the following.

\ct\label{th.acc-cu-gibbs} Let $(\Phi,(\th_\ep)_{\ep>0})$ be
a random perturbation of a $C^2$ diffeomorphism $f$ with an
attractor having a dominated splitting, which is
non-uniformly expanding along the center-unstable direction
and mostly contracting along the center-stable direction.

If $f$ is non-uniformly expanding along the central
direction for random orbits, then every weak$^*$
accumulation point $\mu^0$ of $(\mu^\ep)_{\ep>0}$, when
$\ep\to0$, is a \emph{Gibbs cu-state} for $f$.  \ft

If this holds, then every weak$^*$ accumulation point of the
stationary measures is in the convex hull of the finitely many SRB
measures of $f$.  This means that $f$ is stochastically stable.

In the next sections we prove Theorem~\ref{th.acc-cu-gibbs}.
The strategy is to adapt the notion of \emph{cylinder}
from~\cite{ABV} and~\cite{Vz} to the random setting. A
cylinder $\C$ is the image under a diffeomorphism $\phi$ of
a product of a pair of disks $B^u\times B^s$ of dimensions
$u,s$ respectively, such that the image of any $u$-disk
$\phi\big(B^u\times\{z\}\big), z\in B^s$, is a disk
``almost'' tangent to the center-unstable direction (this is
technically stated using invariant cone fields in the next
sections). The arguments in~\cite{ABV,Vz} (see
also~\cite[Chapter 11]{BDV}) show that in order to obtain a
Gibbs $cu$-state for $f$, it is enough to obtain an
invariant measure $\mu$ with a cylinder $\C$ such that
\begin{enumerate}
\item $\C$ has positive $\mu$-measure;
\item $\mu$ over $\C$ disintegrates along the partition of
  $u$-disks on measures absolutely continuous with respect to
  the volume form on these disks induced by $m$, and
\item the $u$-disks are unstable manifolds for $f$,
  i.e. they are uniformly contracted backwards by $f$.
\end{enumerate}
To prove Theorem~\ref{th.acc-cu-gibbs} we show that in our
setting a stationary ergodic measure $\mu^\ep$ admits a
\emph{cylinder with mass uniformly bounded away from zero}
having $u$-disks \emph{with diameter uniformly bounded away
  from zero} and which are \emph{uniformly contracted
  backward} by a sequence of perturbations.  This is the
content of Section~\ref{sec:centre-unst-cylind}.

Having this, it is not difficult to show that these
cylinders accumulate, when $\ep\to0$, to cylinders having
the same properties with respect to any weak$^*$
accumulation point of $(\mu^\ep)_{\ep>0}$ and for the
unperturbed map $f$. This is the purpose of
Subsections~\ref{sec:absol-cont-limit}
and~\ref{sec:almost-every-ergodic}.

Now we just have to show that every $f$-invariant measure
which is a limit of stationary measures admitting a cylinder
as above must be a Gibbs $cu$-state, which is proved in
Subsection~\ref{sec:absol-cont-accum} following the
arguments in~\cite{ABV,Vz}. Combining these steps we prove
Theorem~\ref{th.acc-cu-gibbs} which together with
Theorems~\ref{t.ABV-prop6.3} and~\ref{t.Vz-ThmA+CorD}
complete the proof of the sufficient condition on
Theorem~\ref{te:mainB}.


\section{Curvature and distortion control}
\label{sec:main-lemmas}


Here we outline some local geometrical and dynamical
consequences of the dominated decomposition and
non-uniformly expanding assumptions, referring mainly to
previous works~\cite{ABV} and~\cite{AA} for proofs.


\subsection{Curvature of center-unstable disks}
\label{sec:bounded-curvature} It was shown in~\cite{ABV}
that $f$ satisfies a bounded curvature property over disks
having the tangent space at each point contained in a cone
field around the center-unstable direction. Here we present
the ``random version'' of the main results in \cite[Section
2]{ABV}.

  Given $0<a<1$ we define the \emph{center-unstable cone
field $C^{cu}_a=(C^{cu}_a(x))_{x\in U}$ of width a} by
\begin{equation}
  \label{eq:cu-cone-field}
  C^{cu}_a(x) = \{
v_1+v_2\in E^{cs}_x\oplus E^{cu}_x \colon \| v_1 \| \le a
\| v_2 \| \}
\end{equation}
and the \emph{center-stable cone field
$C^{cs}_a=(C^{cs}_a(x))_{x\in U}$ of width a} in the same
manner but reversing the roles of the bundles
in~(\ref{eq:cu-cone-field}).

Up to increasing $\la$ slightly we may fix $a$ and $\ep$
small enough so that condition (a) of
Subsection~\ref{sec:part-hyperb-diff} (dominated
decomposition) extends to vectors in the cone fields for
all maps nearby $f$, i.e.
\begin{equation}
  \label{eq:domination-cone-fields}
  \| Df_t(x) v^{cs} \|\cdot\| Df_{t}^{-1}(f_t x) v^{cu}\|
 \le\lambda \|v^{cs}
  \|\cdot \|v^{cu}\|
\end{equation}
for all $v^{cs}\in C^{cs}_a(x)$, $v^{cu}\in C^{cu}_a(f_t
x)$, $x\in U$ and $t\in\supp(\th_\ep)$. Moreover, the
domination property above together with the continuity of
$\Phi$ and the closeness of $t$ to $t^*$ imply that $Df_t
C_a^{cu}(x)$ is contained in a cone of width $\la a$
centered around $Df_t E_x^{cu}$, defined as above with
respect to the splitting $Df_t E_x^{cs} \oplus Df_t
E_x^{cu}$. Since the subspaces $Df_t E_x^{cs}, Df_t
E_x^{cu}$ are close to $E_{f_t x}^{cs}, E_{f_t x}^{cu}$
respectively, then $ Df_t C_a^{cu}(x) \subset C_a^{cu}(f_t
x) $ if $\ep>0$ is small enough. By analogous arguments we
get $D(f_t)^{-1} C_a^{cs}(x)\subset C_a^{cs}(f_t^{-1}x)$.

Given an embedded sub-manifold $S\subset U$ we say that $S$ is
\emph{tangent to the center-unstable cone field} if $T_x S\subset
C^{cu}_a(x)$ for all $x\in S$. Hence $f_t(S)$ is also tangent to
the center-unstable cone field.    The curvature of these
sub-manifolds and their iterates will be approximated in local
coordinates by the notion of H\"older variation of the tangent
bundle as follows.

Let us take $\de_0$ sufficiently small so that if we take
$V_x=B(x,\de_0)$, then the exponential map $\exp_x:V_x\to
T_x M$ is a diffeomorphism onto its image for all $x\in M$.
We are going to identify $V_x$ through the local chart
$\exp_x^{-1}$ with the neighborhood $U_x=\exp_x V_x$ of the
origin in $T_x M$. Identifying $x$ with the origin in $T_x
M$ we get that $E_x^{cu}$ is contained in $C^{cu}_a(y)$ for
all $y\in U_x$, reducing $\de_0$ if needed. Then the
intersection of $E^{cs}_x$ with $C^{cu}_a(y)$ is the zero
vector. So if $x\in S$ then $T_y S$ is the graph of a
linear map $A_x(y):E_x^{cu}\to E_x^{cs}$ for $y\in U_x\cap
S$.

For $C>0$ and $\ze\in (0,1)$ we say that the \emph{tangent bundle
of
  $S$ is $(C,\ze)$-H\"older} if
\begin{equation}
  \label{eq:holder-bundle}
  \| A_x(y) \| \le C \dist_S(x,y)^\ze \quad\mbox{for all}\quad
y\in U_x\cap S \qand x\in U,
\end{equation}
where $\dist_S(x,y)$ is \emph{the distance along $S$}
defined by the length of the shortest smooth curve from $x$
to $y$ inside $S$.

Up to choosing  smaller $a>0$ and $\ep>0$ we may assume
that there are $\la<\la_1<1$ and $0<\ze<1$ such that for
all norm one vectors $v^{cs}\in C^{cs}_a(x), v^{cu}\in
C^{cu}_a(x)$, $x\in U$ it holds $$ \| Df_t(x) v^{cs}
\|\cdot\| Df_{t}^{-1}(f_t x) v^{cu}\|^{1+\ze}
 \le\lambda_1.
$$ For these  values of $\la_1$ and $\ze$, given a $C^1$
sub-manifold $S\subset U$ tangent to the center-unstable
cone field we define
\begin{equation}
  \label{eq:curvature-definition}
  \kappa(S)=\inf\{
C>0 : TS \mbox{ is }(C,\ze)\mbox{-H\"older} \}.
\end{equation}

The proofs of the results that we present below may be
obtained by mimicking the proofs of the corresponding ones
in \cite{ABV}, and we leave them as an exercise to the reader.
The basic ingredients in those proofs are the cone
invariance and dominated decomposition properties that we
have already extended for nearby perturbations $f_t$ of the
diffeomorphism $f$.

\cpr\label{pr:bounded-curvature}
There is $C_1>0$ such that
for every $C^1$ sub-manifold $S\subset  U$ tangent to the
center-unstable cone field and every $\un t\in T^\NN$
\begin{enumerate}
\item there exists $n_1$ such that $\kappa(f^n_{\un t} S)\le
  C_1$ for all $n\ge n_1$ with $f^k_{\un t}S\subset U$ for
  all $1\le k\le n$;
\item if $\kappa(S)\le C_1$ then $\kappa(f^n_{\un t} S)\le
  C_1$ for all $n\ge1$ such that $f^k_{\un t}S\subset U$ for
  all $1\le k\le n$;
\item in particular, if $S$ is as in the previous item, then
$$ J_n: f^n_{\un t} S \ni x \mapsto \log | \det (Df\vert
T_x f^n_{\un t} S) | $$ is $(L_1,\zeta)$-H\"older
continuous with $L_1>0$ depending only on $C_1$ and $f$,
for every $n\ge1$.
\end{enumerate}
\fpr

The bounds provided by
Proposition~\ref{pr:bounded-curvature} may be seen as
bounds on the curvature of embedded disks tangent to the
center-unstable cone field.


\subsection{Hyperbolic times}
\label{sec:bounded-distortion}

>From the condition of non-uniform expansion along the
center-unstable direction we will be able to deduce some uniform
expansion at certain times which are precisely defined through the
following notion.

\cd\label{de:hyperbolic-times} Given $0<\al<1$ we say that
$n\ge1$ is a $\al$-hyperbolic time for $(\un t,x)\in
T^\NN\times U$ if $$ \prod_{j=n-k+1}^n \| Df^{-1} \vert
E^{cu}_{f^j_{\un t} x} \| \le \al^k \quad\mbox{for
all}\quad k=1,\dots,n. $$ \fd

The main technical result ensuring the existence of hyperbolic times
is due to Pliss~\cite{Pli}, whose proof can be found in \cite[Lemma
3.1]{ABV} or \cite[Section 2]{Man}.

\cle \label{l.pliss} Let $H\ge c_2 > c_1$ and
$\zeta={(c_2-c_1)}/{(H-c_1)}$. Given real numbers $a_1,\ldots,a_N$
satisfying $$ \sum_{j=1}^N a_j \ge c_2 N \qand a_j\le H
\;\;\mbox{for all}\;\; 1\le j\le N, $$ there are $l>\zeta N$ and
$1<n_1<\ldots<n_l\le N$ such that $$ \sum_{j=n+1}^{n_i} a_j \ge
c_1\cdot(n_i-n) \;\;\mbox{for
  each}\;\; 0\le n < n_i, \; i=1,\ldots,l. $$
\fle

Using Lemma~\ref{l.pliss} it is not difficult to show that the
condition of non-uniform expansion for random orbits along the
center-unstable direction is enough to ensure that almost all points
have infinitely many hyperbolic times according to the following
result whose proof can be easily adapted from \cite[Corollary
3.2]{ABV}.

\cpr \label{pr:infinite-hyperbolic-times} There exist
$\ga,\al>0$ depending only on $f$ such that for
$\theta^\NN\times\m$ almost all $(\un t,x)\in T^\NN\times U$
and a sufficiently big integer $N\ge1$, there exist $1\le
n_1<\dots<n_k\le N$, with $k\ge\ga N$, which are
$\al$-hyperbolic times for $(\un t,x)$.  \fpr

Let  $n$ be a $\al$-hyperbolic time for $(\un t,x)\in
T^\NN\times U$. This implies  that $Df^{-k}\vert
E^{cu}_{f^n_{\un t} x}$ is a contraction for all
$k=1,\ldots,n$. In addition, if $a>0$ and $\ep>0$ are taken
small enough in the definition of the cone fields and the
random perturbations, then taking $\de_1>0$ also small, we
have by continuity
\begin{equation}
  \label{eq:spare-contraction}
  \| Df_t^{-1} \vert E^{cu}_{f_t y} \| \le \al^{-1/2}
  \| Df^{-1} \vert E^{cu}_{f x} \|
\end{equation}
for all $t\in\supp(\th_\ep^\NN), x\in\ov{f U}$ and $y\in U$ with
$\dist(x,y)<\de_1$. As a consequence of this the next result is
obtained following \cite[Lemma 2.7]{ABV}.

\cle\label{le:contraction-disks} Given any $C^1$ disk
$\De\subset U$ tangent to center-unstable cone field, $x\in
\De$ and $n\ge1$ a $\al$-hyperbolic time for $(\un t,x)$,
we have
$$
 \dist_{f^{n-k}_{\un t}\De}( f^{n-k}_{\un t} (y)  ,
f^{n-k}_{\un t} (x) )\le \al^{k/2} \dist_{f^n_{\un t} \De}(f^n
_{\un t}y, f^n_{\un t} x), \quad k=1,\dots,n,
$$
for every point $y\in \De$ such that
$\dist_{f_{\un t}^n(\De)}
(f^n_{\un t} (y), f^n_{\un t} (x))\le \de_1$.
\fle

Using the previous lemma and the H\"older continuity property given
by Proposition~\ref{pr:bounded-curvature}  the following bounded
distortion result can be deduced as in \cite[Proposition 2.8]{ABV}.

\cpr \label{p.distortion} There exists $C_2>1$ such that,
given any $C^1$ disk $\Delta$ tangent to the center-unstable
cone field with $\kappa(\Delta) \le C_1$, and given any
$x\in \Delta$ and $n \ge 1$ a $\alpha$-hyperbolic time for
$(\un t,x)$, then
$$
\frac{1}{C_2} \le \frac{|\det
  Df^{n}_{\un t} \vert T_y \Delta|} {|\det Df^{n}_{\un t}
  \vert T_x \Delta|} \le C_2
$$
for every $y\in \Delta$ such that $\dist_{f^n_{\un
    t}(\Delta)} ( f^n_{\un t}(y) , f^n_{\un t}(x) )\le \delta_1$.
\fpr



\section{Center-unstable cylinders}
\label{sec:centre-unst-cylind}

Now we show that $\mu^\ep$ admits a cylinder with very
specific properties in the setting of maps with dominated
splitting which are non-uniformly expanding along the
center-unstable direction for random orbits.

Let $\mu^\ep$ be a physical measure of level $\ep$ for some
small $\ep>0$ and take $m_D$ the normalized Lebesgue measure
on some $C^1$ disk $D$ tangent to the center-unstable cone
field such that $m_D$-almost every point of $D$ is in
$B(\mu^\ep)$ and
satisfies \eqref{eq:rNUE}.
It is possible to choose such a disk, because $B(\mu^\ep)$
has nonempty interior in $U$. Now define for each $n\geq 1$
\begin{equation}
  \label{eq:unt}
\mu_n^{\un t} = \frac1n\sum_{j=0}^{n-1} (f^j_{\un t})_{*} \m_D  .
\end{equation}
We know from Theorem~\ref{te:mainA} that each $\mu^\ep$ is
the weak$^*$ limit of the sequence $(\mu_n^{\un t})_n$ for a
$\th_\ep^\NN$ generic $\un t$ by item (2) of
Theorem~\ref{te:mainA}. We fix a $\th_\ep^\NN$ generic $\un
t$ in everything that follows within this section.


A cylinder $\C\subset M$ is the image of a $C^1$
diffeomorphism $\phi: B^u\times B^s\hookrightarrow M$ where $B^k$ is the
$k$-dimensional unit ball of $\RR^k, k=s,u$. We will say
that a $C^1$ disk $D$ \emph{crosses} $\C$ if $D\cap\C$ is a
graph over $B^u$: there exists $g:B^u\to B^s$ such that
$D\cap\C=\{\phi(w,g(w)): w\in B^u\}$.

The following is the main result of this subsection.

\cpr\label{prop.cylinder} Let $\mu^\ep$ be a stationary
probability measure for $\{\Phi,(\th_\ep)_{\ep>0}\}$ where $f$
is a non-uniformly hyperbolic $C^2$ diffeomorphism. Then
there are $\be=\be(f,c_u)>0, \rho=\rho(f,c_u)>0 , d=d(f,c_u)>0$, a cylinder
$\C=\phi(\B^u\times B^s)$ and a family
$\K$ of disks tangent to the center-unstable cone field
which cross  $\C$ and whose union is the set $K$ such that
\begin{enumerate}
\item $\mu^\ep(K\cap\C)\ge \be$ and both $\phi(B^u\times 0)$ and
  $\phi(u\times B^s)$ are disks containing a sub-disk with
  radius $\ge\rho$ for all $u\in B^u$;
\item for every disk $\ga\in\K$ there exists a sequence $\un
  s\in\supp\th_\ep^\NN$ such that $(f_{\un s}^n)^{-1}\mid
  \ga$ is a $\al^{n/2}$-contraction: for $w,z\in\ga$
\[
\dist_{(f_{\un s}^n)^{-1}\ga}\big((f_{\un s}^n)^{-1}(w),(f_{\un
  s}^n)^{-1}(z)\big) \le \al^{n/2}\dist_\ga(w,z)
\]
  where $\dist_\ga$ is the induced distance on $\ga$ by the
  Riemannian metric on $M$.
\item there exists a component $\nu$ of $\mu^\epsilon$ with
  mass uniformly bounded from below by $\beta$ such that
the disintegration $\{\nu_\ga\}_\ga$ of
  $\nu\mid \C$ along the disks $\ga\in\K$ has densities
  with respect to the Lebesgue induced measure $m_\ga$ on
  $\ga$ uniformly bounded from above and below: $d^{-1} \le
  (d\nu_\ga / d m_\ga) \le d$, $\nu_\ga$
  almost everywhere and for almost every $\ga\in\K$.
\end{enumerate}
\fpr

The proof follows an idea in~\cite{ABV}: to consider a
component of the average $\mu_n^{\un t}$ calculated at
hyperbolic times and its weak$^*$ accumulation points.

\begin{proof}
To control the densities of the push-forwards at
hyperbolic times we set
\[
A=\{ x\in D: \dist_D(x,\partial D)\ge\de_1\}
\]
where $\dist_D$ is the distance along $D$, and take $\de_1$
small enough so that $m_D(A)>0$. Then we define for each
$n\ge1$ (we recall that $\un t$ is $\th_\ep^\NN$ generic
fixed from the beginning)
\[
H_n=\{ x\in A: n\mbox{ is a simultaneous hyperbolic time
  for } (\un t,x)\}.
\]
We note that Lemma~\ref{le:contraction-disks} ensures that
$\dist_{f^n_{\un t}(D)}(f^n_{\un t}(x),\partial f^n_{\un
  t}(D))\ge\de_1$ for every $x\in H_n$.

Let $D_n(x,\delta_1)$ be the $\de_1$-neighborhood of
$f^n_{\un t}(x)$ inside $f^n_{\un t}(D)$. Then
Proposition~\ref{p.distortion} ensures that the density of
$((f^n_{\un t})_* m_D)\mid D_n(x,\delta_1)$ with respect to
$m_{D_n(x,\delta_1)}$ is uniformly bounded from above and
from below if we normalize both measures.

To extend this control of the density to a significant portion of
$D$ we use the following result proved in \cite[Proposition 3.3 and
Lemma 3.4]{ABV}.

\cle\label{lem.bolinhas}
There is $\omega>0$ (depending only on $M$, the curvature of
center-unstable disks and on the dimension $u$ of the
center-unstable bundle) such that for all $n\ge1$ we can find a
finite subset $\hat H_n$ of $H_n$ satisfying
\begin{enumerate}
\item $\hat B_n=\{ D_n(x,\de_1/4), x\in \hat H_n\}$ is a pairwise
  disjoint collection;
\item the union $ B_n=\cup \hat B_n$ is such that
  $((f^n_{\un t})_* m_D)( B_n)\ge \omega \cdot m_D(H_n)$.
\end{enumerate}
\fle

\subsection{A special component of the time average}
\label{sec:spec-comp-time}

We define a component of the average measure $\mu_n^{\un t}$ defined
in \eqref{eq:unt}
\begin{equation}
  \label{eq:component}
\nu_n=\frac1n\sum_{j=0}^{n-1} ((f^j_{\un t})_* m_D)\mid  B_j.
\end{equation}

\cle\label{lem:boundedmass} There is $\be_0>0$ such that
$\nu_n(\cup_{j=0}^{n-1} f^j_{\un t}(D))\ge\beta_0$ for all
big enough $n\ge1$.  \fle
\begin{proof}
  We note that
\[
\frac1n\sum_{j=0}^{n-1} m_D(H_j) =\int\int \chi_{H_j}(x) \,
d m_D(x) \, d\#_n(j) = \int\left( \int \chi_{H_j}(x) \,
  d\#_n(j) \right) d m_D(x)
\]
where $\#_n$ is the uniform distribution on
$\{0,\dots,n-1\}$. By
Proposition~\ref{pr:infinite-hyperbolic-times} for big $n$
we must have that the inner integral is bounded from below
by $\ga>0$. By Lemma~\ref{lem.bolinhas} the mass of $\nu_n$
is bounded from below by $\omega \cdot
n^{-1}\sum_{j=0}^{n-1} m_D(H_j)\ge\omega\gamma m_D(D)$ for
big enough $n$. We just have to take $\be_0=\omega\gamma
m_D(D)$ since $\supp(\nu_n)\subset\cup_{j=0}^{n-1} f^j_{\un
  t}(D)$.
\end{proof}

With these settings the support of $\nu_n$ is a finite union
$\cup_{j=0}^{n-1} B_j$ of disks having size bounded from
above and below. Let $\nu$ be an accumulation point of
$(\nu_n)_{n\ge1}$ in the weak$^*$ topology:
$\nu=\lim_k\nu_{n_k}$. Then the support of $\nu$ is
contained in $B_\infty=\cap_{n\ge1}\ov{\cup_{j>n} B_j}$.

This construction shows that for $y\in B_\infty$ there are
sequences $k_j\to\infty$ of integers and disks
$D_j=D_{k_j}(x_{k_j},\de_1/4)$ and points
$y_j\in D_j$ such that $y_j\to y$ when $j\to\infty$.  We
know from subsections~\ref{sec:bounded-curvature}
and~\ref{sec:bounded-distortion} that $D_j$ are $C^1$
center-unstable disks containing a inner $\de_1$-ball.
Moreover the sequence $(D_j)_j$ is relatively compact by the
Ascoli-Arzela Theorem, hence up to taking subsequences we
have $x_{k_j}\to x$ and $D_j\to D_x$ in the $C^1$ topology
when $j\to\infty$ for some $x\in B_\infty$ and a disk $D_x$
centered at $x$ with radius $\de_1/4$.  Thus
$y\in\ov{D_x}\subset B_\infty$.


\subsection{Special sequence of backward contracting parameters}
\label{sec:spec-sequ-param}

Now we obtain the special sequence of parameters for which
we have uniform backward contraction.
Let $(j(n))_{n\ge1}$ be the subsequence of indexes such that
$D_n=D_{k_{j(n)}}\to D_x$ as above when $n\to\infty$. Then
$(t_{k_{j(n)}})_n$ admits a convergent subsequence to some
$s_1\in\supp\th_\ep$.

To avoid too many subscripts we
let that subsequence be indexed by
$k_n^0$ with $n\ge1$. This is a subsequence of $(k_j)_j$. By
definition of hyperbolic times we know that
$(f_{t_{k_n^0}})^{-1}$ is a $\al^{1/2}$-contraction on
$D_{k_n^0}$ for all $n\ge1$.
Hence by the $C^1$
convergence of the disks and the $C^2$ continuity of the
family $\Phi$, we must have that $(f_{s_1})^{-1}$ is a
$\al^{1/2}$-contraction on $D_x$.

We also have that $(t_{k_n^0-1})_n$ admits a subsequence
tending to some $s_2\in\supp\th_\ep$ indexed by $(k_n^1)_n$,
which is a subsequence of $(k_n^0)_n$. In general we have
that $t_{k_n^{\ell-1}-\ell}\to s_\ell$ when $n\to\infty$
where $(k_n^\ell)_n$ is a subsequence of $(k_n^{\ell-1})_n$
for every $\ell\ge0$. The same continuity arguments as above
ensure that $(f_{s_\ell}\circ\dots\circ f_{s_1})^{-1}$ is a
$\al^{\ell/2}$-contraction on $D_x$.

This shows that for every accumulation disk $D_x\in
B_\infty$ as above there exists a subsequence $\un
s\in\supp\th_\ep^\NN$ such that $(f_{\un s}^j)^{-1}\mid D_x$
is a $\al^j$-contraction for every $j\ge1$. We have proved
item (2) in the statement of Proposition~\ref{prop.cylinder}.

\medskip

In what follows, we denote by $\B$ the family of
center-unstable disks in $B_\infty$ obtained through this
limit process.

\subsection{Construction of the cylinder}
\label{sec:constr-cylind}

Now we start the construction of the cylinder. Given any
disk $D\in\B$, the compactness of $B_\infty$ and the
uniformity of $\de_0$ (the radius of invertibility of the
exponential map of $M$ defined in
Subsection~\ref{sec:bounded-curvature}) enables us to
construct a (open) cylinder $\C$ over any sub-disk $D_0$ of
$D$ with radius $\rho\in(0,\min\{\de_0,\de_2\})$ by
considering the images under the exponential map of vectors
in $T_z M$ orthogonal to $T_z D_0$ and with norm less than
$\rho$. We assume that the connected components $\upsilon$
of every center-unstable disk $\ga$ that crosses $\C$ have
diameter smaller than $2\rho$ inside $\ga$. We call $\C$ a
$\rho$-cylinder.

We assume that $\rho<\de_1/100$. Write
$B_j(\rho)$ for the disks obtained from $B_j$
removing the $\rho$-neighborhood of the boundary of every
disk in $B_j$ and let $\hat B_j(\rho)$ be the union of the
points in $B_j(\de)$. Then setting
\[
\nu_{n,\rho}=\nu_n\mid \cup_{j=0}^{n-1} \hat B_j(\de)
\]
we see that $((f^j_{\un t})_* m_D)\mid\hat B_j(\rho) \ge
(1-\de)\cdot ((f^j_{\un t})_* m_D)\mid \hat B_j$ for some
$\de=\de(\rho)>0$.  The value of $\de>0$ may be taken
independently of $j$ because the bounded distortion property
at hyperbolic times (Proposition~\ref{p.distortion}) implies
that the relative mass removed from the disks is comparable
for all iterates.

Hence for a sufficiently small $\rho>0$ as above we may
assume that $\nu_{n,\rho}(M)\ge\be_0/2$ for all $n$ big
enough. We fix this value of $\rho>0$ from now on. Letting
$\nu_\rho$ be an accumulation point of $\nu_{n,\rho}$ for a
subsequence of $(\nu_{n_k})_k$ we have $\nu_\rho\le\nu$.

The $\rho$-cylinders $\C$ constructed as above have uniform
size (depending on $\rho$ only), meaning that they contain a
ball of radius $\rho$.  Observe that
$B_\infty\subseteq\Lambda_\epsilon\subseteq
B(\Lambda,\overline{\epsilon})$ for some
$\overline{\epsilon}>0$, where we write
$B(\Lambda,\overline{\epsilon})$ for
$\cup_{x\in\Lambda}B(x,\overline{\epsilon})$. Then the cover
of $B_\infty$ by the family of all cylinders admits a
minimal cover $\C_1,\dots,\C_k$.

We claim that $k$ is bounded above uniformly independently
of $\epsilon$. Indeed, any cover of $B_\infty$ by such
cylinders is part of a cover of $\Lambda_\epsilon$ by
$\rho$-balls. Let $N$ be the minimum number of
$\rho$-balls needed to cover $\Lambda_\epsilon$. Since
$\Lambda_\epsilon$ is in a $\overline\epsilon$-neighborhood of
$\Lambda$, $N$ is a constant equal to the minimum
number of $\rho$-balls needed to cover $\Lambda$, for small
enough $\epsilon>0$. Hence $k\le N$.

This shows that for some cylinder $\C\in\{\C_1,\dots,\C_k\}$
we must have
\[
\nu(\C)\ge\nu_\rho(\C)\ge
\frac{\nu_\rho(B_\infty)}N=\frac{\nu_\rho(M)}N
\ge \frac{\beta_0}{2N}.
\]
According to the construction of the $\rho$-cylinders, for
every disk $D(\rho)\in\hat B_j(\rho)$ such that $D(\rho)\cap
\C\neq\emptyset$, then the components of $D\cap \C$ cross
$\C$, where $D$ is the corresponding disk in $\hat B_j$
whose truncation gives $D(\rho)$, $j\ge1$. Moreover by an
arbitrarily small change in $\rho$ we may assume that
$\nu(\partial \C)=0$.

Let us denote by $\K_n$ the components of the intersection
$D\cap \C$ that cross $\C$, for all $D\in B_n$ and
$n\ge1$, and let $K_n=\cup\K_n$ be the union of the points
in $\K_n$.  In addition let $\K$ be the set of disks from
$\B$ that cross $\C$ and $K$ the set of all points in $\K$.
Then for all $n\ge1$
\[
\nu_{n,\rho}(\C)=\nu_{n,\rho}(\C\cap\cup_{j=0}^{n-1} K_j)
\le \nu_n(\C\cap\cup_{j=0}^{n-1} K_j).
\]
Hence taking limits of subsequences (recall that
$\nu(\partial \C)=0$ and $\nu_\rho\le\nu$) we arrive at
\[
\frac{\beta_0}{2N}\le \nu_\rho (\C) \le
\limsup_{n\to\infty} \nu_n(\C\cap\cup_{j=0}^{n-1} K_j).
\]
But since $K$ contains the set of accumulation points of
$(\cup_{j=0}^{n-1} K_j)_{n\ge1}$ and $\nu_n$ is defined by
the  average~\eqref{eq:component}, we have that
\[
\limsup_{n\to\infty} \nu_n(\C\cap\cup_{j=0}^{n-1} K_j)
\le \nu(\C\cap K)
\]
and so $\nu(\C\cap K)\ge \beta_0/(2N)$.

However $\mu^\ep\ge\nu$ by construction, hence
$\mu^\ep(\C\cap K)\ge \beta_0/(2N)$ also. We stress that
either $\beta_0, N$ or $\rho$ do not depend on the choice of
$\ep$ nor of $\un t$.

This proves the statement of item (1) of
Proposition~\ref{prop.cylinder}.


\subsection{Densities along center-unstable disks}
\label{sec:dens-along-centre}

Let $\C$ be the $\rho$-cylinder constructed before and let
$D_0\in\K$ be the base disk of $\C$. Write $\K_j$ and
$K_j$ as before and set $\K_\infty=\K$. We take a sequence
$(\P_n)_{n\ge1}$ of increasing partitions of the family
$\D=\cup_{0\le k\le\infty} K_j$ as follows.

The cylinder $\C$ is endowed with the orthogonal projection
onto the base disk $p:\C\to D_0$ and the disks in $\D$
define a projection $\pi:\D\to\xi_0$ where $\xi_0=
p^{-1}(\{x_0\})$ for some fixed $x_0\in D_0$.

To define $\P_k$, first take a sequence $(\V_k)_{k\ge1}$ of
increasing partitions of $\xi_0$ whose diameter tends to
zero. Next introduce the space $\hat\D=\cup_{0\le
  j\le\infty} K_j\times\{j\}$.  Then for any given $k$ we
say that two elements $(x,i)\in K_i\times\{i\}$ and
$(y,j)\in K_j\times\{j\}$ of $\hat\D$ are in the same atom
of $\P_k$ if both $x,y$ project under $\pi$ into the same
atom of $\V_k$ and either $i,j\ge k$ or $i=j<k$.

Observe that since $\xi_0$ is diffeomorphic to a ball of
some Euclidean space and we may identify each disk
$\gamma\in\K_j$ with $\pi(\gamma)$ for all $0\le
j\le\infty$, then we may assume that the union $\partial
\P_k$ of the boundaries of the elements of $\P_k$ satisfies
$\mu(\partial\P_k)=0$ by suitably choosing the sequence
$\V_k$, i.e., the boundaries of the elements of $\V_k$
should have zero measure with respect to
$\hat\mu=\pi_*(\mu)$.

Given $x\in\hat\D$ and writing $\P_k(x)$ for the atom of
$\P_k$ which contains $x$, it is clear that
$\P_k(x)\supset\P_{k+1}(x)$ for all $k\ge1$ and also that
$\cap_{k\ge1}\P_k(x)$ equals $\pi^{-1}(\{x\})\cap\C$.

Let $A$ be a Borel subset of $D_0$ and
$\zeta\in\K_j$. Since the projection
$p$ sends $\zeta$ diffeomorphically on $D_0$ and in $\C$
the angles involved in the projection are uniformly bounded,
we may find a constant $C>0$ such that
\begin{equation}
  \label{eq:dens1}
  \frac1{C}\cdot\frac{m_0(A)}{m_0(D_0)}
  \le
  \frac{m_{\zeta}\big( p^{-1}(A)\cap \zeta
    \big)}{m_{\zeta}(\zeta)}
  \le
  C \cdot\frac{m_0(A)}{m_0(D_0)},
\end{equation}
where we have written $m_0$ and $m_{\zeta}$ for the
Lebesgue induced measures on $D_0$ and $\zeta$ by $m$,
respectively.

Proposition~\ref{p.distortion} ensures that the density
of $\big(f^j_{\un t}\big)_*(m_D)$ with respect to Lebesgue
measure on each disk $\gamma\in\K_j$ is bounded from
above and from below, thus
\begin{equation}
  \label{eq:dens2}
  \frac1{C_2}
  \le
  \frac{\big(f^j_{\un t}\big)_*(m_D)\big( p^{-1}(A)\cap
    \zeta \big)}{m_{\zeta}\big( p^{-1}(A)\cap \zeta
    \big)}
  \le
  C_2
  \quad\mbox{and}\quad
  \frac1{C_2}
  \le
  \frac{\big(f^j_{\un t}\big)_*(m_D)\big(
    \zeta \big)}{m_{\zeta}(\zeta)}
  \le
  C_2.
\end{equation}
Combining (\ref{eq:dens1})
and (\ref{eq:dens2}) we get
\begin{equation}
  \label{eq:dens3}
  \frac1{C_2^2 C}\cdot\frac{m_0(A)}{m_0(D_0)}
  \le
  \frac{\big(f^j_{\un t}\big)_*(m_D)\big( p^{-1}(A)\cap \zeta
    \big)}{\big(f^j_{\un t}\big)_*(m_D)(\zeta)}
  \le
  C_2^2 C \cdot\frac{m_0(A)}{m_0(D_0)},
\end{equation}
for all big enough values of $j$.

Now we define a sequence $(\hat\nu_n)_{n\ge1}$ of measures
on $\hat D$ by
\[
\hat\nu_n\big(E_0\times\{0\}\cup \dots\cup
E_{n-1}\times\{n-1\}\big)
=\frac1n\sum_{j=0}^{n-1} (f_{\un t}^j)_*(m_D)(E_j)
\]
where $E_i\subset K_i$ for $i=0,\dots,n-1$,
and $\hat\nu_n(E)=0$ for all $E\subset\cup_{n\le j\le\infty}
K_j$. Observe that given $k\ge1$ and $x\in\hat\D$ the atom
$\P_k(x)$ is formed by a union of disks in $\cup_{0\le j\le
  k-1}\K_j$. Hence by the definition of $\hat\nu_n$ and
by~(\ref{eq:dens3}) we conclude that
\begin{equation}
  \label{eq:dens4}
  \frac1{d}\cdot\hat\nu_n(\P_k(x))\cdot m_0(A)
  \le
  \hat\nu_n\big(p^{-1}(A)\cap\P_k(x)\big)
  \le
  d\cdot\hat\nu_n(\P_k(x))\cdot m_0(A),
\end{equation}
where $d=C_2^2 C / m_0(D_0)$.

It is easy to see that any weak$^*$ accumulation point of
$\hat\nu_n$ is supported in $K\times\{\infty\}$.
Moreover if we choose a sequence $n_k$ such that
$\nu_{n_k}\to\nu$, then this just means that $\hat\nu_{n_k}$
tends to a measure $\hat\nu$ such that
$\hat\nu(E\times\{\infty\})=\nu(E)$ for all $E\subset K$.
Since we may assume without loss that
$\nu\Big(\partial\big(\P_k(D)\cap p^{-1}(A)\big)\Big)=0$ for all
$k\ge0$, by the choice of $\V_k$ during the construction
above, the inequalities~(\ref{eq:dens4}) also hold in the
limit, i.e.
\begin{equation}
  \label{eq:dens5}
  \frac1{d}\cdot\hat\nu(\P_k(x))\cdot m_0(A)
  \le
  \hat\nu\big(p^{-1}(A)\cap\P_k(x)\big)
  \le
  d\cdot\hat\nu(\P_k(x))\cdot m_0(A).
\end{equation}
By the theorem of Radon-Nikodym (\ref{eq:dens5}) means that
the density of $\nu$ along the disks $\gamma\in\K$ is
bounded above and below, as stated in item (3) of
Proposition~\ref{prop.cylinder}.  This concludes the proof
of this proposition.
\end{proof}


\section{Accumulation cylinders}
\label{sec:accum-cylind}

In what follows we fix a decreasing sequence $\ep_k\to0$
when $k\to\infty$ and a sequence $\mu_k=\mu^{\ep_k}$ of
ergodic stationary measures. We write also $\nu_k$ for the
component of each $\mu_k$ for $k\ge1$ with well behaved
disintegrations given by Proposition~\ref{prop.cylinder}.

We observe that since $(\supp\th_\ep)_{\ep>0}$ is a nested
family of connected compact subsets shrinking to $\{t^*\}$
when $\ep\to0$, and for each $\ep>0$ and any stationary
measure $\mu^\ep$ the set $\supp\mu^\ep$ is $f_t$-invariant
for all $t\in\supp\th_\ep$, we may choose the sequence
$\mu_k$ so that $(\supp\mu_k)_k$ is a nested family of
$f$-invariant compact subsets.

\cpr\label{prop.accylinder}
Let $\mu$ be a weak$^*$ accumulation point of
$(\mu_k)_k$. Then there exists $d_0>0$, a cylinder $\C$, a family
$\K$ of disks tangent to the center-unstable cone field
which cross  $\C$, whose union is the set $K$, and a
component $\nu$ of $\mu$ such that
\begin{enumerate}
\item $\nu(K\cap\C)\ge \be$;
\item  $(f^n)^{-1}\mid \ga$ is a $\al^{n/2}$-contraction on
  every disk $\ga\in\K$, where $f=f_{t^*}$.
\item the disintegration $\{\nu_\ga\}_\ga$ of $\nu\mid \C$
  along the disks $\ga\in\K$ has densities with respect to
  the Lebesgue induced measure $m_\ga$ on $\ga$ uniformly
  bounded from above and below: $d_0^{-1} \le (d\nu_\ga / d
  m_\ga) \le d_0$; 
\end{enumerate}
\fpr

The value of $\beta$ above is the same from the
statement of Proposition~\ref{prop.cylinder}. The value of
$d_0$ depends only on $d$ from Proposition~\ref{prop.cylinder}.

Item (2) above shows that every  $\gamma\in\K$ is a
center-unstable disk in $U$. This means that $E^{cu}_x$ is
uniformly expanded by $Df$ for every $x\in\gamma$. The
domination property for the splitting $E^{cu}\oplus E^{cs}$
guarantees that any eventual expansion along the
complementary direction is weaker than this. Thus $\gamma$ is
contained in the unique local strong-unstable manifold
$W^u_{loc}(x)$ tangent to $E^{cu}_x$, see~\cite{Pe76}.

In item (3) we can also say that the disintegration of $\mu$
along the disks of $\K$ has densitities bouded from above
and from below, since $\nu$ is a component of $\mu$.

\begin{proof}[Proof of Proposition~\ref{prop.accylinder}]
  Let $\mu_k$ be as stated in the beginning of the
  subsection and let $\C_k, \K_n$ and $K_n$ be the corresponding
  cylinders, families of disks and sets from
  Proposition~\ref{prop.cylinder}. We assume that
  $\mu_k\to\mu$ in the weak$^*$ topology when $k\to\infty$.
  Then $\mu$ is an $f$-invariant probability measure
  (Remark~\ref{re.accinvariant}). We take also $\nu$ a limit
  point of $\nu_k$ in the weak$^*$ topology. Note that since
  $\nu_n\le\mu_n$ for all $n\ge1$ then $\nu\le\mu$ also.

  The compactness of $M$ ensures that for some subsequence
  $k_n$ the cylinder $\C_{k_n}$ tends to a cylinder $\C$. In
  fact, each cylinder $\C_k$ is a diffeomorphic image of
  $\phi_k:B^u\times B^s \hookrightarrow M$, with $B^\ell$
  the $\ell$-dimensional unit ball of $\RR^\ell, \ell=s,u$.
  By the Ascoli-Arzela Theorem there is a subsequence
  $(k_n)_{n\ge1}$ such that $\phi_{k_n}(B^u\times 0)$
  converges in the $C^1$-topology to a disk
  $D_0=\phi(B^u\times 0)$ in $M$. Since the diameters of
  $\phi_{k_n}(B^u\times 0)$ and $\phi_{k_n}(0\times B^s)$
  are uniformly bounded from below by $\rho>0$ (by
  Proposition~\ref{prop.cylinder}) and by the construction
  of $\C_k$, defining $\C$ as the set of images under the
  exponential map of vectors in $T_z M$ orthogonal to $T_z
  D_0$ and with norm less than $\rho$, then $\ov{\C}_{k_n}$
  tends to $\ov{\C}$ in the Hausdorff topology.

Let $\K$ be the family of disks $D$ in $\C$ which are
accumulated by sequences of disks $D_n$ in $\K_{k_n}$ for
$n\ge1$. Since every disk $D_n$ is tangent to the
center-unstable cone field of $f$, the continuity of the
cone field on $U$ assures that every disk $D\in\K$ is also a
center-unstable disk.

\cre\label{rmk:crossing}
It will be useful to note that up to taking a slightly
smaller base disk $D_0$ we may assume without loss that the
disks in $K_{k_n}$ cross $\C$ for all big enough $k$.
\fre

For any fixed $\gamma\in \K$ let $x,y\in D$ and take
$(x_n)_n,(y_n)_n$ sequences in $\gamma_n\in\K_{k_n}$ such that
$x_n\to x$ and $y_n\to y$ when $n\to\infty$. From item (3)
of Proposition~\ref{prop.cylinder} we know that there are
sequences of parameters $(\un s(n))_{n\ge1}$ such that $\un
s(n)\in\supp\th_{\ep_{k_n}}^\NN$ and
\[
\dist_{(f_{\un s(n)}^j)^{-1}(\gamma_n) }\big((f_{\un
  s(n)}^j)^{-1}(x_n),(f_{\un s(n)}^j)^{-1}(y_n)\big) \le
\al^{j/2}\dist_{\gamma_n}(x_n,y_n)
\]
for every $j\ge1$ and for every given $n\ge1$. Fixing
$j\ge1$  we get
\[
(s_1(n),\dots,s_j(n)) \to (t^*,\dots,t^*)\quad
\mbox{when}\quad n\to\infty,
\]
because $\supp(\th_{\ep_{k_n}}^\NN)\to\{t^*\}$.
The continuity of $f_t(x)$ with respect to $(t,x)\in T\times
M$ implies that
\begin{equation*}
\dist_{f^{-j}(\gamma) }\big(f^{-j}(x),f^{-j}(y)\big) \le
\al^{j/2}\dist_{\gamma}(x,y)
\end{equation*}
for every given $j\ge1$. Hence $f^{-j}$ is an
$\al^{j/2}$-contraction on every $\gamma\in\K$, which proves item
(2) of Proposition~\ref{prop.accylinder}.

Let now $K$ be the union of the elements of $\K$. By
construction we have the following property
\begin{align}\label{eq:property0}
  \text{for all $\de>0$ there is $n_0\in\NN$ such that
    $K_{k_n}\subset B(K,\de)$ for all $n\ge n_0$. }
\end{align}
Since $\nu_k(\C_k\cap\K_k)\ge\beta>0$ for all $k\ge1$
from Proposition~\ref{prop.cylinder},
letting $\de>0$ be such that $\mu(\partial B(K,\de))=0$ and
so also $\nu(\partial B(K,\de))=0$
(this holds except for an at most countable set of values of
$\de$), then
\[
\nu(B(K,\de))=\lim_{n\to\infty} \nu_{k_n}(B(K,\de)) \ge \be >0.
\]
Moreover $K=\cap_{\de>0} B(K,\de)$ thus $\nu(K)=\inf_{\de>0}
\nu(B(K,\de))\ge\beta$.  This proves item (1) of the
statement of Proposition~\ref{prop.accylinder}.


\subsection{Absolute continuity of limit  measure on
  accumulation cylinder}
\label{sec:absol-cont-limit}

Here we prove item (3) of Proposition \ref{prop.accylinder}.
We recall that the limit cylinder $\C$ has base $D_0\in\K$.

We take a sequence $(\P_n)_{n\ge1}$ of increasing partitions
of the family $\hat\D\subset\cup_{1\le j\le\infty}
K_j\times\{j\}$ of all disks which cross $\C$ (by
Remark~\ref{rmk:crossing} this family contains disks from
infinitely many distinct $\K_k$), defined in the same
fashion as in the proof of
Proposition~\ref{prop.cylinder}. For the rest of the proof
we write $\K_\infty$ for $\K$.

Let $p:\C\to D_0$ be the orthogonal projection on the base
disk and $\pi:\D\to\xi_0$ be the projection along the leaves
of $\cup_j\K_j$ where $\xi_0= p^{-1}(\{x_0\})$ for
some fixed $x_0\in D_0$. Take a sequence
$(\V_k)_{k\ge1}$ of increasing partitions of $\xi_0$ with
diameter tending to zero and define $\P_k$ in the same way
as in Subsection~\ref{sec:dens-along-centre}.

Exactly as in Subsection~\ref{sec:dens-along-centre}, we may
assume without loss that the union $\partial \P_k$ of the
boundaries of the elements of $\P_k$ have zero measure with
respect to $\hat\mu=\pi_*(\mu)$, by an adequate choice of
the sequence of partitions. Note that given $x\in\D$, then
$\P_k(x)\supset\P_{k+1}(x)$ for all $k\ge1$ and
$\cap_{k\ge1}\P_k(x)$ equals $\pi^{-1}(\{x\})\cap\C$.

Let $A$ be a Borel subset of $D_0$ and $\zeta\in\K_j$. Then
we have (\ref{eq:dens1}) by the same reasons.  Moreover item
(3) of Proposition~\ref{prop.cylinder} ensures that the density
of $\nu_j$ with respect to Lebesgue measure on $\zeta$ is
bounded from above and from below, thus
\begin{equation}
  \label{eq:proj2}
  \frac1{d}
  \le
  \frac{\nu_j\big( p^{-1}(A)\cap \zeta \big)}
  {m_{\zeta}\big( p^{-1}(A)\cap \zeta \big)}
  \le
  d
  \quad\mbox{and}\quad
  \frac1{d}
  \le
  \frac{\nu_j(\zeta_i)}{m_{\zeta}(\zeta)}
  \le
  d.
\end{equation}
Combining (\ref{eq:dens1}) with (\ref{eq:proj2}) we get for
all $\zeta\in\K_j$ and all $j$
\begin{equation}
  \label{eq:proj3}
  \frac1{ C d^2}\cdot\frac{m_0(A)}{m_0(D_0)}
  \le
  \frac{\nu_j\big( \pi^{-1}(A)\cap \zeta
    \big)}{\nu_j(\zeta)}
  \le
  C d^2 \cdot\frac{m_0(A)}{m_0(D_0)}.
\end{equation}
Likewise the argument in
Subsection~\ref{sec:dens-along-centre} we define a sequence
$(\hat\nu_n)_{n\ge1}$ of measures on $\hat\D$ by
\[
\hat\nu_n\big(E_0\times\{0\}\cup \dots\cup
E_{n-1}\times\{n-1\}\big)
=\frac1n\sum_{j=0}^{n-1} \nu_j(E_j)
\]
where $E_i\subset K_i$ for $i=0,\dots,n-1$,
and $\hat\nu_n(E)=0$ for all $E\subset\cup_{n\le j\le\infty}
K_j$. Note that for $k\ge1$ and $x\in\hat\D$ the atom
$\P_k(x)$ is a union of disks in $\cup_{0\le j\le
  k-1}\K_j$. Hence by the definition of $\hat\nu_n$ and
by~(\ref{eq:proj3}) we deduce
\begin{equation}
  \label{eq:proj4}
  \frac1{d_0}\cdot\hat\nu_n(\P_k(x))\cdot m_0(A)
  \le
  \hat\nu_n\big(p^{-1}(A)\cap\P_k(x)\big)
  \le
  d_0\cdot\hat\nu_n(\P_k(x))\cdot m_0(A),
\end{equation}
where $d_0= C d^2/ m_0(D_0)$. Note that $D_0$ is an
accumulation disk and so its size depends only on the value
of $d$, hence $d_0$ depends only on $d$.

As in the proof of item (3) of
Proposition~\ref{prop.cylinder}, any weak$^*$ accumulation
point of $\hat\nu_n$ is supported in $K\times\{\infty\}$
and $\nu_{n}\to\nu$  means that
$\hat\nu_{n}$ tends to a measure $\hat\nu$ such that
$\hat\nu(E\times\{\infty\})=\nu(E)$ for all $E\subset K$.
Then the inequalities~(\ref{eq:proj4}) also hold in the
limit
\begin{equation*}
  \frac1{d_0}\cdot\hat\nu(\P_k(x))\cdot m_0(A)
  \le
  \hat\nu\big(p^{-1}(A)\cap\P_k(x)\big)
  \le
  d_0\cdot\hat\nu(\P_k(x))\cdot m_0(A),
\end{equation*}
which, by the theorem of Radon-Nikodym means that
the density of $\nu$ along the disks $\gamma\in\K$ is
bounded above and below, since $A$, $k$ and $x$ are arbitrary.

This proves item (3) of Proposition~\ref{prop.accylinder}
and concludes the proof.
\end{proof}


\subsection{Almost every ergodic component is a Gibbs
   {\em cu}-state}
\label{sec:almost-every-ergodic}

Here and in the next subsection we conclude the proof of
Theorem~\ref{th.acc-cu-gibbs} by showing that \emph{every
  probability measure~$\mu$, obtained as a limit measure of
  the measures $\mu_k$ (which is $f$-invariant, see
  Remark~\ref{re.accinvariant}), and admitting a component
  $\nu$ satisfying properties (1)-(3) of
  Proposition~\ref{th.acc-cu-gibbs} must be a Gibbs
  $cu$-state.}

Let $(\mu_x)_{x\in M}$ be the ergodic decomposition of $\mu$
(see e.g. \cite{Man}), that is
\begin{itemize}
\item $\mu_x$ is an $f$-invariant \emph{ergodic} probability
  measure for $\mu$-a.e. $x\in M$, and
\item $\int\vfi\,d\mu=\int\!\big(\int
  \vfi\,d\mu_x\big)\,d\mu(x)$ for every
  bounded measurable function $\vfi:M\to\RR$.
\end{itemize}
This family is uniquely defined up to a subset of $M$ with
zero $\mu$-measure.

Let now $\mu$ be a weak$^*$ accumulation point of the
sequence $(\mu_k)_{k\ge1}$ as in the previous subsection.

\cle
\label{le.qtpGibbs}
The measure $\mu_x$ is a  $cu$-Gibbs state for
$\mu$-a.e. $x\in K$.
\fle

\begin{proof}
  Let $\Sigma$ be the full $\mu$ measure subset of $M$ where
  the family $(\mu_x)$ giving the ergodic decomposition of
  $\mu$ is defined. Let $R$ be the subset of \emph{Oseledets
  regular points} with respect to $f$ (see e.g. \cite{Man}),
that is, the full $\mu$-measure subset of $M$ on whose
points Lyapunov exponents are well defined. Then
$K\cap\Sigma\cap R$ has full $\mu$-measure on $K$, since
$\mu(K)>0$ (recall the statement of
Proposition~\ref{th.acc-cu-gibbs}).

Observe that by the uniform backward contraction property of
the disks of $\K$ we have for $\mu$-a.e. $x\in
K\cap\Sigma\cap R$
\[
\lim_{n\to+\infty}\frac1n \log\| Df^{-n}\mid E^{cu}_x\| \le
\frac12\log\alpha < 0,
\]
which implies that the Lyapunov exponents along the
tangent direction $E^{cu}$ to the disks of $\K$ are all positive.

Consider a measurable set $B_0$ such that
\[
m_\gamma(B_0\cap\gamma)=0\quad\mbox{for all}\quad \gamma\in\K
\]
and $\mu(B_0)$ is maximal among all measurable subsets with
this property. Since $\nu$ is absolutely continuous along
the disks of $\K$ we have $\nu(B_0)=0$. In what follows set
$Z_0=K\cap\Sigma\cap R\setminus B_0$.

By definition of ergodic decomposition, for every measurable
set $A\subset Z_0$ we have
\begin{equation}
  \label{eq:ergdecomp1}
\mu(A)=\int\mu_x(A)\,d\mu(x)
\end{equation}
and we want to express this as an integral over $Z_0$ in
order to use the properties of $\nu$.

The ergodicity of $\mu_x$ ensures that $
\mu_x(A)=\lim_{n\to+\infty} n^{-1}\sum_{j=0}^{n-1}
\chi_A(f^j(x)) $ for $\mu$-a.e. $x$, where $\chi_A(x)=1$ if
$x\in A$ and $0$ otherwise.  Hence except for a $\mu$-zero
subset of points we see that $\mu_x(A)>0$ only if $x$ has
some iterate in $A\subset Z_0$.

Let $k(z)$ be the smallest positive integer such that
$f^{-k(z)}(z)\in Z_0$, which is defined $\mu$-almost
everywhere in $Z_0$. Note that $\mu_z=\mu_{f^i(z)}$ for all
$i\in\ZZ$.  Thus we may write
\[
\mu(A)=\int_{Z_0} k(z)\,\mu_z(A)\,d\mu(z),
\]
since for $\mu$-a.e. $z$ we remove
$\mu_z(A)=\mu_{f^{-1}(z)}(A)=\dots=\mu_{f^{-k(z)+1}(z)}(A)$
from (\ref{eq:ergdecomp1}).

Now we use the following technical result which can be deduced from
§3 of \cite{Rk62}, and whose proof
  can also be found in \cite[Lemma 6.2]{ABV}.

\cle
\label{l.conditionalmeasures}
Let $\lambda$ be a finite measure on a measure space $Z$,
with $\lambda(Z)>0$.  Let $\K$ be a measurable partition of
$Z$, and $(\lambda_z)_{z\in Z}$ be a family of finite
measures on $Z$ such that
\begin{enumerate}
\item the function $z\mapsto \lambda_z(A)$ is measurable,
  and it is constant on each element of $\K$, for any
  measurable set $A\subset Z$
\item $\{w:\lambda_z = \lambda_w\}$ is a measurable set with
  full $\lambda_z$-measure, for every $z\in Z$.
\end{enumerate}
Assume that $\lambda(A) = \int \ell(z) \lambda_z(A) \,d\lambda$ for
some measurable function $\ell:Z\to\RR_+$ and any measurable subset
$A$ of $Z$.  Let $\{\tilde\lambda_\gamma$, $\gamma\in\K\}$, and
$\{\tilde\lambda_{z,\gamma}$, $\gamma\in\K\}$, be disintegrations of
$\lambda$ and $\lambda_{z}$, respectively, into conditional
probability measures along the elements of the partition $\K$.  Then
\( \tilde\lambda_{z,\gamma} = \tilde\lambda_{\gamma} \) for
$\lambda$-almost every $z\in Z$ and $\hat\lambda_z$-almost every
$\gamma$, where $\hat\lambda_z$ is the quotient measure induced by
$\lambda_z$ on $\K$. \fle

We set $Z=Z_0$, $\lambda=(\mu\mid Z_0)$, $\lambda_z=(\mu_z\mid
Z_0)$, $\ell(z)=k(z)$ and $\K$ as before, with $z\in Z_0$, and apply
Lemma~\ref{l.conditionalmeasures}. We conclude that the
disintegration $\mu_{z,\gamma}$ of $\mu_z$ along the disks
$\gamma\in\K$ coincides almost everywhere with the disintegration
$\mu_\gamma$ of $\mu\mid Z_0$ along the same family of disks.
Therefore, since we have already shown by
Proposition~\ref{prop.accylinder} and the choice of $B_0$ that
$\mu_\gamma$ is absolutely continuous with respect to $m_\gamma$,
and also that the Lyapunov exponents of $\mu_z$ along the tangent
directions to the disks of $\K$ are all positive, we conclude that
$\mu_z$ is an ergodic $cu$-Gibbs state for $f$, for $z\in Z_0$. This
finishes the proof of the lemma.
\end{proof}



\subsection{The accumulation measure is a Gibbs \emph{cu}-state}
\label{sec:absol-cont-accum}

Here we finish the proof of Theorem~\ref{te:mainB}. We fix
$\ep_k\to0$, $\mu_k=\mu^{\ep_k}$ and $\mu=\lim_{k\to\infty}\mu_k$ in
the weak$^*$ topology as in the previous section. We denote by $\C$
the cylinder and by  $K$ the compact subset in the statement of
Proposition~\ref{prop.accylinder} with respect to $\mu$. We also
denote by $\C^k$ the cylinder and by $K^k$ the compact subset from
the statement of Proposition~\ref{prop.cylinder} with respect to
each $\mu_k$, $k\ge1$.

Define $G$ to be the set of all points $x\in\Sigma\cap
R$ (i.e. Oseledets regular points whose orbit defines an
ergodic $f$-invariant measure, see the proof of
Lemma~\ref{le.qtpGibbs}) such that $\mu_x$ is a $cu$-Gibbs
state and set
\[
\nu_0=\int_G \mu_x \, d\mu(x).
\]
Since $\mu_x=\mu_{f^i(x)}$ for all $i\in\ZZ$ and $\mu$-a.e.
$x$ the measure $\nu_0$ is $f$-invariant and not identically
zero since $G\supset K$. By construction, $\nu_0/\nu_0(G)$
is a $cu$-Gibbs state.

The purpose of this section is to prove the next result.

\cpr\label{prop.abscontmu} $\mu=\nu_0$, that is $\mu$ is a
Gibbs $cu$-state.  \fpr

Before the proof we present some useful lemmas.  By the
results of \cite{BV,ABV} and by Theorem~\ref{t.ABV-prop6.3}
we know that $\nu_0/\nu_0(G)$ is a physical probability
measure and we can use the following result which
corresponds to \cite[Lemma 3.5]{Vz}. For each $n\ge1$ set
$K_n=\{ x\in M: \tau(x)\le n\}$, where
$\tau(x)=\min\{k\ge1:f^{-k}(x)\in K\}$.  Note that by
Poincar\'e's Recurrence Theorem the function $\tau$ is
finite almost everywhere with respect to any $f$-invariant
probability measure which assigns positive mass to $K$.

\cle\label{le.serieuniforme}
Let $R_j=\{x\in K: \tau(x)=j\}$.
Then there exists $C>0$ and $\lambda_0\in(0,1)$ such that
given any physical measure $\nu$ for $f$ we have
$\nu(R_j) \le C \cdot \lambda_0^j$
for every $j\ge1$.
\fle

Now for every $k,n\ge1$ we set
\[
A_n^k=\bigcap_{\un t\in\supp(\th_{\ep_k}^\NN)}
\bigcap_{j=1}^{n} f_{\un t}^j\big(M\setminus K^k\big)\cap K^k
\]
and define $K_n^k=M\setminus A_n^k$ and $\mu_k^n=\mu_k\mid K_n^k$.

The next result is a consequence of the assumption that $f_t$ is is
$C^2$-close to $f\equiv f_{t^*}$ when $t$ is close to $t^*$ together
with Lemma~\ref{le.serieuniforme}.

\cle\label{lem:randomtail}
There is $C_0>0$ and for any given $n\in\NN$ there is $\ell\ge1$ such that
\[
\mu_k(A_n^k)<C_0\lambda_0^{n+1}
\quad\mbox{for every}\quad k>\ell.
\]
Moreover for every fixed $n\ge1$ we have
$\mu_k^n\to\mu\mid K_n$ when $k\to\infty$ in the weak$^*$ topology.
\fle

\begin{proof}
  We start by noting that by Lemma~\ref{le.qtpGibbs} we have
  $\mu\mid K = \nu_0\mid K$. So from
  Lemma~\ref{le.serieuniforme} for any given $n$ we have
  that
  \[
  A_n=\bigcap_{j=1}^n f^j(M\setminus K) \cap K
  \]
  is such that $x\in A_n$ if, and only if, $f^{-j}(x)\in
  M\setminus K$ for $j=1,\dots,n$ and $x\in K$, and this is
  equivalent to $\tau(x)>n$, that is $A_n\subseteq\cup_{k>n}
  R_k$.  Thus $\mu(A_n) \le C\sum_{j>n}
  \lambda_0^j=C^\prime\lambda_0^{n+1}$, where
  $C^\prime=C/(1-\lambda_0)$.

  Using the property \eqref{eq:property0} from the
  construction of $K$, since $K$ is closed and $\mu$ is a
  Borel regular measure, fixing the number $n$ of iterates
  involved we have that for small enough $\zeta_1,\zeta_2>0$
  the set
  \[
  A_n(\zeta_1,\zeta_2)=\bigcap_{j=1}^n f^j\big( B(M\setminus
  K, \zeta_2) \big) \cap B(K,\zeta_1)
\]
also satisfies $\mu(A_n(\zeta_1,\zeta_2))< 2
C^\prime\lambda_0^{n+1}$, where $B(X,\zeta)=\cup_{x\in X}
B(x,\zeta)$ is the $\zeta$-neighborhood of $X$ in $M$ for
any $\zeta>0$ and any subset $X$.  Moreover
$A_n(\zeta_1,\zeta_2)$ is an open neighborhood of $A_n$ so
through an arbitrarily small change in $\zeta_1,\zeta_2$ we may
assume that $\mu(\partial A_n(\zeta_1,\zeta_2) )=0$ in what
follows.

Note that for fixed $n$ and $\zeta_1$ as above, there exists
$\ell_1\in\NN$ such that for small enough $\zeta_2<\zeta_1$
we have $K^k\subset B(K,\zeta_2/2)\subset B(K,\zeta_1)$ and
$M\setminus K^k\subset B(M\setminus K,\zeta_2)$ for all
$k>\ell_1$ (recall \eqref{eq:property0}).  Since $f_t$
depends continuously on $t$ and the number $n$ of iterates
is fixed, we may take $\ell_1$ big enough and $\zeta_2>0$
small enough such that
$A_n^k \subset A_n(\zeta_1,\zeta_2)$
for all $k\ge\ell_1$.

But $\mu_k\to\mu$ in the weak$^*$ topology so by the
assumption on the boundary of $A_n(\zeta_1,\zeta_2)$ we
arrive at $\mu_k(A_n^k)\le \mu_k(A_n(\zeta_1,\zeta_2))\le 4
C^\prime\lambda_0^{n+1}$ for all big enough $k$. The first
statement of the lemma is obtained resetting $\ell$ to a
bigger value (if needed) and letting $C_0=4 C^\prime$.

For small enough $\xi_1,\xi_2>0$ we now define
\[
A_n^k(\xi_1,\xi_2)= \bigcap_{\un t\in\supp(\th_{\ep_k}^\NN)}
\bigcap_{j=1}^{n} f_{\un t}^j\big(B(M\setminus K^k,
\xi_2)\big)\cap B(K^k,\xi_1)
\]
which is an open neighborhood of $A_n^k$. Again by the
construction of $K$ and by \eqref{eq:property0}, for fixed
$n$ and $\xi_1>0$ there exists $\ell_2\in\NN$ such that for
small enough $\xi_2<\xi_1$ we have $K\subset
B(K^k,\xi_2/2)\subset B(K^k,\xi_1)$ for all $k>\ell_2$.
Since $f_t$ is $C^2$ close to $f=f_{t^*}$ for big $k$ we
can check that $A_n\subset A_n^k(\xi_1,\xi_2)$ for all
$k\ge\ell_2$.

Choosing small values of $\ze_1,\ze_2>0$ and taking $\ell$
big enough the above arguments ensure that for all
$k\ge\ell$ it holds
\begin{equation}
  \label{eq:encaixe}
A_n\subset A_n(\zeta_1,\zeta_2)\subset
A_n^k(2\ze_1,2\ze_2)\subset A_n(4\zeta_1,4\zeta_2)
\end{equation}
and we can simultaneously assume that we have
\begin{equation}
  \label{eq:fronteira0}
  \mu(\partial A_n(\zeta_1,\zeta_2) )= 0 = \mu(\partial
  A_n(4\zeta_1,4\zeta_2)).
\end{equation}
Let us take an open set $B$ such that $\mu(\partial B)=0$
(the collection of all such sets generates the Borel
$\sigma$-algebra $\mu\bmod0$). Using \eqref{eq:encaixe} we
get for all $k\ge\ell$
\[
\mu_k(B\cap [M\setminus A_n(4\zeta_1,4\zeta_2)])
\le
\mu_k(B\cap [M\setminus A_n^k(2\ze_1,2\ze_2)])
\le
\mu_k(B\cap [M\setminus A_n(\zeta_1,\zeta_2)])
\]
and letting $k\to\infty$ and using~\eqref{eq:fronteira0} we
arrive at
\begin{eqnarray*}
  \mu(B\cap [M\setminus A_n(4\zeta_1,4\zeta_2)])
& \le &
\liminf_{k\to\infty} \mu_k(B\cap [M\setminus
A_n^k(2\ze_1,2\ze_2)])
\\
& \le &
\limsup_{k\to\infty} \mu_k(B\cap [M\setminus
A_n^k(2\ze_1,2\ze_2)])
\\
&\le&
\mu(B\cap [M\setminus A_n(\zeta_1,\zeta_2)]).
\end{eqnarray*}
Finally letting $\ze_2\to0$ first and then $\ze_1\to0$ also
the compact sets $M\setminus A_n(4\zeta_1,4\zeta_2)$ and
$M\setminus A_n(\zeta_1,\zeta_2)$ both grow to $M\setminus
A_n$ which clearly equals $K_n$. In the same way $M\setminus
A_n^k(2\ze_1,2\ze_2)$ grows to $M\setminus A_n^k=K_n^k$.

This together with the last sequence of inequalities shows
that $\mu_k^n(B)\to(\mu\mid K_n)(B)$ when $k\to\infty$,
finishing the proof.
\end{proof}

  For each $n\ge1$ set $\nu_n=\nu_0\mid K_n$.  Note that
  $\nu_n=\mu\mid K_n$ since $\nu_0\mid K=\mu\mid K$ and both
  $\mu$ and $\nu_0$ are $f$-invariant.  In the weak$^*$
  topology $\nu_n\to\nu_0$ when $n\to\infty$ because of the
  following simple fact.

\cle\label{le.returns}
For $\nu_0$-almost every $z$ there exists $j \le 0$ such that
$f^j(z)\in K$.
\fle

\dem Given $B$ a Borel subset we have that
$\mu_x(B)=\lim_{n\to+\infty}n^{-1}\sum_{j=0}^{n-1}
\chi_B(f^j(x))$ for $\mu$-almost every $x$. Since
$\mu(K)>0$, if $\nu_0(B)>0$ then for some $x\in K$ we have
$\mu_x(B)>0$ and there exists $j\ge0$ such that $f^j(x)\in
B$.  \cqd

Note that Lemma~\ref{lem:randomtail} implies that
$\mu^n_k\to\mu_k$ when $n\to\infty$ in the weak$^*$ topology
in a uniform way, since $\mu_k(A_n^k)\to0$ when $n\to\infty$
uniformly in $k$.  In addition, Lemma~\ref{lem:randomtail}
also ensures that for any fixed $n\in\NN$ we have
$\mu^n_k\to\nu_n$ when $k\to\infty$ in the weak$^*$ topology
by the definition of $\nu_n$. We are now ready to complete
the proof of Proposition~\ref{prop.abscontmu}.

\begin{proof}[Proof of Proposition~\ref{prop.abscontmu}]
Let $\zeta>0$ and a continuous $\vfi:M\to\RR$ be given. Then
we may find a big enough $n\ge1$ such that
\[
|\nu_0(\vfi)-\nu_n(\vfi)| \le \zeta\quad\mbox{and}\quad
|\mu^n_k(\vfi)-\mu_k(\vfi)| \le\zeta
\]
for every sufficiently big $k\ge1$. Observe that for
the right hand inequality above we need the uniform bound on
$\mu_k(A_n^k)$ provided by Lemma~\ref{lem:randomtail}.

Having fixed $n$ we may now take $k$ big enough keeping the
above inequalities and satisfying also
\[
|\mu_k(\vfi)-\mu(\vfi)| \le\zeta
\quad\mbox{and}\quad
|\nu_n(\vfi)-\mu^n_k(\vfi)|\le\zeta.
\]
Finally putting this all together we arrive at
\begin{eqnarray*}
  |\nu_0(\vfi)-\mu(\vfi)|
  &\le&
  |\nu_0(\vfi)-\nu_n(\vfi)|+|\nu_n(\vfi)-\mu^n_k(\vfi)|+
\\
& &
  |\mu^n_k(\vfi)-\mu_k(\vfi)|+|\mu_k(\vfi)-\mu(\vfi)|
 \le  4\zeta.
\end{eqnarray*}
This concludes the proof of
Proposition~\ref{prop.abscontmu}.
\end{proof}

As explained in the beginning of
Subsection~\ref{sec:sufficient} this is precisely what is
needed to conclude stochastic stability for $f$.
Theorem~\ref{te:mainB} is proved.

\section{A stochastically stable class}
\label{sec:stoch-stable-class}

In this section we present a robust class of partially hyperbolic
diffeomorphisms satisfying conditions (a)-(c) and also
condition~\eqref{eq:rNUE} for random orbits. Here we take $U$
equal to $M$. This presentation follows closely~\cite{ABV} and we
just sketch the main points. The $C^1$ open classes of transitive
non-Anosov diffeomorphisms presented in \cite[Section~6]{BV}, as
well as other robust examples from \cite{Man78}, are constructed
in a similar way.

We start with a linear Anosov diffeomorphism $\hat f$ on the
$d$-dimension\-al torus $M=\TT^d$, $d\ge 2$. We write $TM=E^u\oplus
E^s$ the corresponding hyperbolic decomposition of the tangent fiber
bundle. Let $V$ be a small closed domain in $M$ for which there
exist unit open cubes $K^0$ and $K^1$ in $\RR^d$ such that $V
\subset \pi(K^0)$ and $\hat f(V)\subset \pi(K^1)$, where
$\pi:\RR^d\to \TT^d$ is the canonical projection. Now, let $f$ be a
diffeomorphism on $\TT^d$ such that
 \begin{enumerate}
 \item[(A)] $f$ admits
invariant cone fields $C^{cu}$ and $C^{cs}$, with small width $a>0$
and containing, respectively, the unstable bundle $E^u$ and the
stable bundle $E^s$ of the Anosov diffeomorphism $\hat f$;
 \item[(B)] $f$ is \emph{volume hyperbolic}:
 there is $\sigma_1>1$ so that $$ |\det(Df\vert
T_x\De^{cu})| > \sigma_1 \quad\mbox{and}\quad |\det(Df\vert
T_x\De^{cs})| < \sigma_1^{-1} $$ for any $x\in M$ and any disks
$\De^{cu}$, $\De^{cs}$ tangent to $C^{cu}$, $C^{cs}$, respectively.
 \item[(C)] $f$ is $C^1$-close to $\hat f$ in
the complement of $V$, so that there exists $\sigma_2<1$ satisfying
 $$
 \|(Df \vert T_x \De^{cu})^{-1}\| < \sigma_2
\quad\mbox{and}\quad \|Df \vert T_x \De^{cs}\| < \sigma_2
 $$
for any $x\in (M\setminus V)$ and any disks $\De^{cu}$, $\De^{cs}$
tangent to $C^{cu}$, $C^{cs}$, respectively.
 \item[(D)] there exist some small $\delta_0>0$ satisfying
 $$
 \|(Df \vert T_x\De^{cu})^{-1}\| <
1+\delta_0\quad\mbox{and}\quad \|Df \vert T_x \De^{cs}\| <
1+\delta_0
 $$
  for any $x\in V$ and any disks $\De^{cu}$ and $\De^{cs}$ tangent to
$C^{cu}$ and $C^{cs}$, respectively.
\end{enumerate}

Closeness in (C) should be enough to ensure that $f(V)$ is
also contained in the projection of a unit open cube.  If
$\tilde f$ is a torus diffeomorphism satisfying (A), (B),
(D), and coinciding with $\hat f$ outside $V$, then any map
$f$ in a $C^1$ neighborhood of $\tilde f$ satisfies all the
previous conditions. Results in~\cite[Appendix]{ABV} show in
particular that for any $f$ satisfying (A)--(D) there exist
$c_u>0$ such that $f$ is non-uniformly expanding along its
center-unstable direction as in condition (b) at
Subsection~\ref{sec:part-hyperb-diff}. Below we reobtain
this result in Lemma~\ref{pr.nonunifexp.localdiffeo} as a
particular case of the proof of non-uniform expansion for
random orbits.  In addition from~\cite{T} we know that these
diffeomorphisms are also mostly contracting, as in condition
(c) at Subsection~\ref{sec:part-hyperb-diff}, and we deduce
this as a particular case of the arguments for random orbits
in what follows.

\ct\label{te:mainC} Let $\D$ be the class of $C^1$
diffeomorphisms $f$ satisfying conditions (A)--(D) above.
Every $f\in\D$ is non-uniformly expanding for random orbits
along the center-unstable direction.  \ft

We prove this result in the next section. Perturbing the
original Anosov diffeomorphism $\hat f$ on $V$ through a
flip bifurcation of a fixed point we can obtain $f$ with a
hyperbolic fixed point with a (stable) index different from
the dimension of sub-bundle $E^{s}$, or even with a source or
a sink with eigenvalues close to one in absolute value.
This shows that the class $\D$ contains open sets of maps
satisfying the assumptions of Theorem~\ref{te:mainB} with
\emph{neither} uniform contraction \emph{nor} uniform
expansion along the sub-bundles of their dominated splitting.


\subsection{Behavior over random orbits}
\label{sec:central-expansion}

Here we take $f$ as described above and $T$ a small neighborhood
of $f$ in the $C^2$ topology in such a way that conditions
(A)--(D) hold for every $g\in T$. We define $\Phi(t)=f_t=t$ and
take $(\th_\ep)_{\ep>0}$ a family of measures on $T$ as before. We
are going to show that any such $f$ is non-uniformly expanding
along the center-unstable direction on random orbits and, in the
process, we will also see that this is a neighborhood of maps with
mostly contracting center-stable direction. This will be done by
showing that there is $c>0$ such that for any disk $\De^{cu}$
tangent to the center-unstable cone field and $\th_\ep^\NN\times
m$ almost all $(\un t,x)\in T^\NN\times\De^{cu}$
\begin{equation}
\limsup_{n\to\infty}\frac{1}{n}\sum_{j=0}^{n-1} \log\|\big(Df \vert
T_{f^j_{\un t}x}\De^{cu}_j(\un t)\big)^{-1}\|\le -c,
 \end{equation}
where $\De^{cu}_j(\un t)=f^j_{\un t}\De^{cu}$. To explain this, let
$B_1, \ldots, B_p, B_{p+1}=V$ be any partition of $\TT^d$ into small
domains, in the same sense as before: there exist open unit cubes
$K_i^0$ and $K_i^1$ in $\RR^d$ such that
\begin{equation}
B_i \subset \pi(K_i^0) \quad\text{and}\quad f(B_i) \subset
\pi(K_i^1).
\end{equation}

\noindent Let us fix $\De^{cu}$  any disk tangent to the
center-unstable cone field and define $m$ to be Lebesgue measure in
$\De^{cu}$ normalized so that $m(\De^{cu})=1$.

\cle \label{l.frequency} There is $\zeta>0$ such that for
$\th_\ep^\NN\times m$ almost all $(\un t,x)\in T^\NN\times \De^{cu}$
and large enough $n\ge1$ we have
\begin{equation} \label{eq.frequency}
\#\{0\le j < n: f^j_{\un t}(x) \in B_1\cup\ldots\cup B_p \} \ge
\zeta  n.
\end{equation}
Moreover, there is $0<\tau<1$ for which the set $I_n$ of points
$(\un t,x)\in T^\NN\times \De^{cu}$ whose orbits do not spend a
fraction $\zeta$ of the time in $B_1\cup\ldots\cup B_p$ up to
iterate $n$ is such that $(\th_\ep^\NN\times m)(I_n)\le \tau^n$ for
large $n\ge1$.

In particular, if we take a constant $\un t=(f,f,f,\dots)$
we get the same conclusion for the unperturbed $f$, that is,
there is $\tau\in(0,1)$ such that the set $J_n$ of points
$x\in \De^{cu}$ whose orbits do not spend a
fraction $\zeta$ of the time in $B_1\cup\ldots\cup B_p$ up to
iterate $n$ is such that $m(J_n)\le \tau^n$ for
large $n\ge1$.
\fle

\dem Let us fix $n\ge1$ and $\un t\in T^\NN$. For a sequence $\un
i=(i_0,\ldots,i_{n-1})\in\{1,\ldots,p+1\}^n$ we write
\[
 [\un i]=\De^{cu}\cap
B_{i_0}\cap (f^1_{\un t})^{-1}(B_{i_1})\cap\cdots \cap (f^{n-1}_{\un
t})^{-1}(B_{i_{n-1}})
\]
and define $g(\un i)=\#\{ 0\le j < n : i_j \le p \}$. We
start by observing that for $\zeta>0$ the number of
sequences $\un i$ such that $g(\un i)<\zeta n$ is bounded by
$$
\sum_{k<\zeta n} {n \choose k} p^k \le \sum_{k\le \zeta
  n} {n \choose k} p^{\zeta n} . $$
Using Stirling's formula
(cf.~\cite[Section 6.3]{BV}) the expression on the right
hand side is bounded by $(e^\gamma p^\zeta)^n$, where
$\gamma>0$ depends only on $\zeta$ and $\gamma(\zeta)\to 0$
when $\zeta\to0$.

Assumption (B) ensures that $m([\un i])\le \si_1^{-(1-\zeta)n}$
(recall that $m(\De^{cu})=1$). Hence the measure of the union $I_n(\un t)$
of all the sets $[\un i]$ with $g(\un i)<\zeta n$ is bounded
by
\[
\si_1^{-(1-\zeta)n} (e^\gamma p^\zeta )^n.
\]
Since $\si_1>1$ we may
choose $\zeta$ so small that $e^\gamma p^\zeta  < \si_1^{(1-\zeta)}$.
Then $m(I_n(\un t))\le \tau^n$ with $\tau=e^{\gamma+\zeta-1} p^\zeta
<1$ for big enough $n\ge N$. Note that $\tau$ and $N$ do not depend
on $\un t$.

\cre\label{rmk.uniformeasure} If $x\in M\setminus I_n(\un t)$, then
\[
\#\{0\le j<n : f^j_{\un t}(x)\in B_{p+1}\}\le (1-\zeta)\cdot n
\]
by definition of $I_n(\un t)$. \fre

Setting
\begin{equation}
  \label{eq:I_n}
I_n=\mcup_{\un t\in T^\NN} \big(\{\un t\} \times
I_n(\un t)\big)
\quad\mbox{we have}\quad (\th^\NN_\ep\times m)(I_n)\le
\tau^n
\end{equation}
for big $n\ge N$, by Fubini's Theorem. Since $\sum_n
(\th^\NN_\ep\times m)(I_n) < \infty$ then Borel-Cantelli's
Lemma implies $$
(\th^\NN_\ep\times m) \left( \mcap_{n\ge1}
  \mcup_{k\ge n} I_k \right) =0 $$
and this means that
$\th^\NN_\ep\times m$ almost every $(\un t,x)\in T^\NN\times
\De^{cu}$ satisfies~(\ref{eq.frequency}).

For the unperturbed case just observe that the estimate in
(\ref{eq:I_n}) is uniform in $\un t$, so it also holds for
a constant $\un t=(f,f,f,\dots)$ and the rest of the
arguments are unchanged.
\cqd

\cle \label{pr.nonunifexp.localdiffeo} For $0<\lambda<1$ there are
$\eta>0$ and $c_0>0$ such that, if $f_t$ also satisfies conditions
(C) and (D) for all $t\in T$, then we have
\begin{enumerate}
\item $m\big(\{ x\in M: \sum_{j=0}^{n-1}
\log\|\big(Df \vert T_{f^j_{\un t}x}\De^{cu}_j(\un t)\big)^{-1}\|\le
-c n \}\big) \ge 1-\tau^n$ for all $\un t\in T^\NN$ and for every
large $n$;
\item $\limsup_{n\to\infty}\frac{1}{n}\sum_{j=0}^{n-1}
\log\|\big(Df \vert T_{f^j_{\un t}x}\De^{cu}_j(\un t)\big)^{-1}\|\le
-c$
for $\th^\NN_\ep\times m$ almost all $(\un t,x)\in T^\NN\times
\Delta^{cu}$.
\end{enumerate}


In particular, for a constant $\un t=(f,f,f,\dots)$ item (2)
above still holds for this $\un t$ and $m$-almost every
$x\in\De^{cu}$.

 \fle

\dem Let $\{ B_1,\ldots,B_p,B_{p+1},\ldots,B_{p+1}\}$ be a
measurable cover of $M$ as before and $\zeta>0$ be the constant
provided by Lemma~\ref{l.frequency}. We fix $\eta>0$ sufficiently
small so that $\lambda^\zeta(1+\eta)\le e^{-c}$ holds for some $c>0$
and take $x\in M\setminus I_n(\un t)$ for some $n\ge1$ and $\un t\in
T^\NN$. Conditions (C) and (D) now imply
\begin{equation} 
\prod_{j=0}^{n-1}\|\big(Df \vert T_{f^j_{\un t}x}\De^{cu}_j(\un
t)\big)^{-1}\| \le \lambda^{\zeta n} (1+\eta)^{(1-\zeta)n} \le e^{-c
n}.
\end{equation}
by Remark~\ref{rmk.uniformeasure}. Hence the set in item 1 is
contained in $M\setminus I_n(\un t)$, proving the statement of this
item by the second part of the statement of Lemma~\ref{l.frequency}.

This also means that the second item of the statement holds for
$\th^\NN\times m$ almost every $(\un t,x)\in T^\NN\times \De^{cu}$ by the
statement of Lemma~\ref{l.frequency}.


For the unperturbed case the arguments are analogous.
\cqd

Since $\Delta^{cu}$ was an arbitrary disk tangent to the
center-unstable cone, we conclude from
Lemma~\ref{pr.nonunifexp.localdiffeo} that $f$ is
non-uniformly expanding along the center-unstable direction
for random orbits.  Moreover the same arguments apply
verbatim to $\|Df^{-1}\mid T_{f^{-j} x} f^{-j}(\De^{cs})\|$
and any disk $\De^{cs}$ tangent to the center-stable cone,
so that Lemmas~\ref{l.frequency}
and~\ref{pr.nonunifexp.localdiffeo} holds for this cone
also.

Since the statements of both  Lemmas~\ref{l.frequency}
and~\ref{pr.nonunifexp.localdiffeo} are true in the
unperturbed case, this shows that every map $f\in T$ is
non-uniformly expanding along the center-unstable direction
and mostly contracting along the center-stable direction.





\end{document}